%% file: MC_revised_RRR.tex
\documentclass[onefignum,onetabnum]{siamart250211}

\usepackage{mathtools}
\usepackage{graphicx}
\usepackage{subcaption}
\newtheorem{assumption}{Assumption}

\newcommand{\by}[1]{{\color{black}#1}}
\newcommand{\byy}[1]{{\color{black}#1}}

\input{ex_shared}

\ifpdf
\hypersetup{
  pdftitle={Convergence Analysis of a Stochastic Interacting Particle-Field Algorithm for 3D Parabolic-Parabolic Keller-Segel System},
  pdfauthor={Boyi Hu, Zhongjian Wang, Jack Xin, Zhiwen Zhang},
  colorlinks=true,
  allcolors=black
}
\fi

\begin{document}

\maketitle

\begin{abstract}
Chemotaxis models describe the movement of organisms in response to chemical gradients. In this paper, we propose a stochastic interacting particle-field algorithm integrated with random batch approximation (SIPF-$r$) for the three-dimensional (3D) parabolic-parabolic Keller-Segel (KS) system, also referred to as the fully parabolic KS system. The SIPF-$r$ method approximates the KS system by coupling particle-based density representations with a smooth field variable computed via spectral methods. By incorporating the random batch method (RBM), we bypass the mean-field limit and significantly reduce computational complexity. Under mild assumptions on the regularity of the original KS system and the boundedness of numerical approximations, we prove that the empirical measure associated with the SIPF-$r$ particle system converges in the $1$-Wasserstein distance to the exact measure of the limiting McKean-Vlasov process with high probability. Finally, we present numerical experiments to validate the theoretical convergence results, and demonstrate the efficiency and robustness of the SIPF-$r$ method as a diagnostic tool for density concentration and potential finite-time singularity in 3D parabolic-parabolic KS systems. \by{Specifically, the SIPF-$r$ method identifies the critical threshold of initial mass leading to blow-up under specific initial data configurations.}
\end{abstract}



\section{Introduction}
Chemotaxis is a biological phenomenon involving the movement of organisms (e.g., bacteria) in response to signals (e.g., chemoattractants), which can be produced by the organisms themselves. Theoretical and mathematical modeling of this phenomenon was initiated by Patlak \cite{patlak1953random} and by Keller and Segel \cite{keller1970initiation}. 
In this work, we focus on the fully parabolic KS system as follows:
\begin{equation}
\begin{split}
    \rho_{t}  &= \nabla \cdot (\mu \, \nabla \rho - \chi\, \rho\,  \nabla c ), \\
    \epsilon \,   c_t  &= \Delta\, c - \lambda^2 \, c + \rho, \qquad
    \mathbf{x} \in \Omega \subseteq \mathbb{R}^d, \quad t \in [0, T],
\end{split}
\label{K-S system}
\end{equation}
where $\chi, \mu $ are positive constants, and $\epsilon, \lambda$ are non-negative constants. The model is called parabolic-elliptic if $\epsilon =0$, and fully parabolic if $\epsilon >0$. Here, $\rho$ denotes the density of active particles (bacteria), and $c$ represents the concentration of chemoattractants emitted by these bacteria. The model finds broad applications in biology, ecology, and medicine \cite{perthame2006transport,painter2019mathematical,bubba2019chemotaxis,SHERRATT1994129}.

The nonlinear, potentially singular behavior of KS equations, especially blow-up phenomena \cite{nagai1995blow, herrero1998self, bellomo2015toward, fatkullin2013study}, makes numerical methods essential. Mesh-based methods such as finite difference \cite{chertock2018high,saito2005notes,epshteyn2012upwind}, finite element \cite{saito2007conservative,epshteyn2009fully,saito2011error}, and finite volume schemes \cite{filbet2006finite,chertock2008second,zhou2017finite} are widely employed. Recent advances include local discontinuous Galerkin methods with optimal convergence rates \cite{li2017local} and semi-discrete symmetrization-based schemes that avoid nonlinear solvers \cite{liu2018positivity}. Despite their success, challenges remain in ensuring stability, convergence, and the effective handling of singularities, making numerical studies on the KS model an active and evolving field of research. 

In addition to mesh-based methods, particle-based approaches have also been developed to address the challenges posed by the KS system, providing a complementary perspective. Stevens \cite{stevens2000derivation} developed an $N$-particle system and established its convergence for the fully parabolic case. Ha{\v{s}}kovec and Schmeiser \cite{havskovec2009stochastic} proposed a convergent regularized particle system for the 2D parabolic-elliptic KS model. Moreover, Godinho and Quininao \cite{godinho2015propagation} showed well-posedness results and the propagation of chaos property in a subcritical KS equation. Craig and Bertozzi \cite{craig2016blob} proved the convergence of a blob method for the related aggregation equation. Liu and Yang \cite{liu2017random} introduced a random particle blob method with a mollified kernel for the parabolic-elliptic case, proving its convergence when the macroscopic mean-field equation possesses a global weak solution \cite{liu2019propagation}. \by{Due to the potential finite-time blow-up of the KS system, deriving \emph{a priori} error estimates inherently requires sufficient regularity of the exact solution. Consequently, theoretical convergence analyses in the literature are typically restricted to the pre-blow-up regime \cite{stevens2000derivation, li2017local, liu2017random, zhou2017finite}. For the singular regime, numerical studies primarily focus on structural preservation (e.g., positivity and energy dissipation) and empirical robustness \cite{chertock2008second, epshteyn2009discontinuous, shen2020unconditionally}.}

In \cite{wang2025novel}, we proposed a novel stochastic interacting particle-field (SIPF) algorithm for the fully parabolic KS system \eqref{K-S system} in 3D. The SIPF method approximates the KS solution $\rho$ as the empirical measure of particles (see Eq.~\eqref{rho_particle}) coupled with a smoother field variable $c$ computed using the spectral method (see Eq.~\eqref{trig_ser}). The algorithm employs an implicit Euler discretization and a one-step recursion based on Green’s function, which efficiently captures concentration behavior and potentially finite-time blow-up with only dozens of Fourier modes in each dimension.

Despite the algorithm's numerical efficiency and stability as observed in \cite{wang2025novel}, a complete convergence analysis remains open. This paper fills that gap by establishing and validating convergence estimates for the SIPF-$r$ method, a random-batch variant designed to reduce computational cost. \by{Our main result, presented in Theorem \ref{theorem}, establishes the convergence of the solution $(\widetilde{\rho}, \widetilde{c})$ obtained from the SIPF-$r$ method to the exact solution $(\rho, c)$ in the pre-blow-up regime under mild assumptions}. Specifically, the $1$-Wasserstein distance between the SIPF-$r$ and exact density distributions, denoted as $\mathcal{W}_1(\widetilde{\rho}_{t}, \rho_{t})$, depends on the time step $\delta t$, the number of Fourier modes $H$, and the number of particles $P$, and scales as $\mathcal{O}\left(\frac{1}{H^2} + \sqrt{\frac{1}{P}} + \sqrt{\delta t}\right)$. The proof, detailed in Section~\ref{section: proof}, builds upon key lemmas that quantify single-step update errors and analyze their propagation over time. \by{Furthermore, we demonstrate that the spectral truncation inherently regularizes the singular kernel, ensuring the unconditional stability of the algorithm even as the exact solution approaches the finite-time singularity.}

The error between the chemical concentration $\widetilde{c}$ in the SIPF-$r$ method and the exact solution $c$ originates from Fourier truncation and implicit Euler time discretization. Under appropriate regularity assumptions on $c$, the truncation error decays as the number of Fourier modes $H$ increases. The temporal error can be expressed through differences between Fourier coefficients of $\widetilde{c}$ and $c$, which further couples with the $L^2$ particle trajectory error $\mathbb{E}(\|\widetilde{X}_t - X_t\|_{L^2})$ accumulated from the preceding step. Here, $X_t$ and $\widetilde{X}_t$ are the exact dynamics of the system associated with the density and its numerical approximations obtained via our methods, respectively. 

At the numerical discretization level, the RBM \cite{jin2020random, jin2021random, huang2025mean, cai2024convergence} is incorporated into the SIPF-$r$ method, which renders the particles fully independent and identically distributed (i.i.d.), thereby effectively circumventing the need to address the propagation of chaos \cite{liu2019propagation}. At each time step, small random batches of particles are selected with replacement for particle interactions. In the error estimate between $\widetilde{c}$ and $c$, the i.i.d. property, together with the complex-valued mean value theorem \cite{markushevich2013theory} and Bernstein's inequality, bound the deviations of the empirical mean of particle trajectories from its expectation. This deviation accounts for the uncertainty described in Theorem \ref{theorem}. The numerical experiments in Section \ref{section: experiment} further demonstrate that, with the introduction of the RBM, the numerical examples maintain a high level of precision.

The error between $X_t$ and $\widetilde{X}_t$ is influenced by the gradient of $c$ and $\widetilde{c}$, reflecting the dependence of the particle trajectories on the interaction potential. Using Parseval's identity \cite{korner2022fourier}, we relate $\|\nabla\widetilde{c}-\nabla c\|_{L^2}$ to Fourier-coefficient errors, establishing two coupled recursive inequalities between $\|\nabla\widetilde{c}-\nabla c\|_{L^2}$ and $\mathbb{E}(\|\widetilde{X}_{t} - X_{t}\|_{L^2})$. The first inequality relates $\|\nabla \widetilde{c}-\nabla c\|_{L^2}$ to $\mathbb{E}(\|\widetilde{X}_t - X_t\|_{L^2})$; the second updates $\mathbb{E}(\|\widetilde{X}_t - X_t\|_{L^2})$ by incorporating a $\|\nabla \widetilde{c}-\nabla c\|_{L^2}$ term (see Eqs.~\eqref{eq:a_n}-\eqref{eq:b_n}). Substituting and decoupling these aforementioned recursive inequalities yields a bound for $\mathbb{E}(\|\widetilde{X}_t - X_t\|_{L^2})$ that depends only on earlier errors. Through the natural coupling $\gamma_t = \text{Law}(\widetilde{X}_t,X_t)$, this $L^2$ error estimate translates to the $1$-Wasserstein distance $\mathcal{W}_1(\widetilde{\rho}_t,\rho_t)$. Applying the discrete Gronwall inequality together with Fourier-coefficient error estimates leads to the global convergence result in Theorem \ref{theorem}.

In addition to the theoretical analysis, Section~\ref{section: experiment} provides comprehensive numerical experiments to validate the performance of the SIPF-$r$ method. These experiments verify the predicted convergence rates, confirm the validity of the regularity assumptions, and explore the complex dynamics of the KS system. \by{A particular focus is placed on identifying the initial-data-dependent threshold for finite-time blow-up in three dimensions under specific initial data configurations.} In Subsection~\ref{subsection:blowup-detection}, the SIPF-$r$ algorithm identifies blow-up regimes and captures severe concentration toward finite-time singularities (see Figures~\ref{fig:blowup-comparison} and~\ref{fig:two-spheres-blowup}). Notably, the method remains robust in detecting blow-up even under practical discretization parameters ($H$, $P$, $\delta t$) that are much coarser than those required for the convergence analysis (see Figure~\ref{fig:c_max_vs_time}).

The rest of the paper is organized as follows. In Section \ref{algorithm}, we review the SIPF method for solving the fully parabolic KS system and present the derivation of the SIPF-$r$ method, which incorporates the RBM to compute the particle interaction. In Section \ref{section: proof}, under certain assumptions, we provide a detailed convergence analysis of the SIPF-$r$ method, which decomposes the proof into several lemmas. In Section \ref{section: experiment}, we present numerical results to validate the necessity of the assumptions, demonstrate the accuracy of the SIPF-$r$ method, and confirm the theoretical convergence rate derived in our analysis. \byy{Section~\ref{section: discussion} discusses extensions of the SIPF-$r$ framework to a wider class of mathematical-biological chemotaxis models, and Section~\ref{section: conclusion} concludes the paper.}

\section{Derivation of the SIPF-$r$ Method}\label{algorithm}
In this section, we present the SIPF-$r$ method for solving the fully parabolic KS model. Since the evolutionary phenomena of interest occur in the interior of the domain, we restrict the system \eqref{K-S system} to a large domain $\Omega=[-L/2,L/2]^3$ and assume Dirichlet boundary conditions for particle density $\rho$ and Neumann boundary conditions for chemical concentration $c$.

Throughout this section, we use the standard notation $\rho$, $c$, etc., to represent the exact solutions of the fully parabolic KS model. For the variables computed or approximated using the SIPF-$r$ algorithm, we instead use the notations $\widetilde{\rho}$, $\widetilde{c}$, etc.

As a numerical algorithm, we assume that the temporal domain $[0,T]$ is partitioned by $\{t_n\}_{n=0:n_T}$ with $t_0=0$ and $t_{n_T}=T$. We approximate the density $\widetilde{\rho}$ at $t=t_n$ by empirical particles $\{\widetilde{X}^{p}_{t_n}\}_{p=1:P}$, i.e.,
\begin{align}
	\widetilde{\rho}_{t_n} = {\frac {M_0}{P}}\, \sum_{p=1}^{P} \delta (x - \widetilde{X}^{p}_{t_n}),  \; \; P\gg 1,
	\label{rho_particle}
\end{align}
where $M_0$ denotes the conserved total mass (integral of $\rho$ over the domain $\Omega$). For the chemical concentration $\widetilde{c}$, we adopt a Fourier basis approximation. Its unnormalized Fourier coefficients are defined by
\[
\widetilde{\alpha}_{t;\mathbf{j}}
:=\mathcal{F}_{\mathbf{j}}[\widetilde c(\cdot,t)]
=\int_{\Omega}\widetilde c(\mathbf{x},t)e^{-i\omega_{\mathbf{j}}\cdot\mathbf{x}}\,d\mathbf{x},
\qquad \omega_{\mathbf{j}}=\frac{2\pi}{L}\mathbf{j}.
\]
The corresponding inverse Fourier representation is
\begin{equation}
	\widetilde{c}(\mathbf{x},t)=\frac{1}{L^3}\sum_{\mathbf{j} \in \mathcal{H}}\, \widetilde{\alpha}_{t;\mathbf{j}}\, \exp(i 2\pi j_1 \, x_1/L) \exp(i 2\pi j_2\, x_2/L)\exp(i2\pi j_3\, x_3/L), \label{trig_ser}
\end{equation}
where $\mathcal{H}$ denotes the index set, i.e.,
\by{\begin{align}\label{index set H}
    \mathcal{H}=\{\mathbf{j}\in \mathbb{Z}^3: |j_1|,|j_2|,|j_3|\leq \frac{H}{2}\},
\end{align}}
and $i=\sqrt{-1}$. Defining $\alpha_{t;\mathbf{j}}:=\mathcal{F}_{\mathbf{j}}[c(\cdot,t)]$ in the same way, the exact solution $c(\mathbf{x},t)$ can also be approximated by a truncated spatial Fourier series expansion as follows:
\begin{equation}\label{c fourier}
c(\mathbf{x},t) \approx \frac{1}{L^3}\sum_{\mathbf{j}\in \mathcal{H}}\, \alpha_{t;\mathbf{j}}\, \exp(i 2\pi j_1 \, x_1/L) \exp(i 2\pi j_2\, x_2/L)\exp(i2\pi j_3\, x_3/L).
\end{equation}

\begin{remark}
The choice of the Fourier basis over Hermite polynomials for approximating the chemical concentration is based on the fact that, since the blow-up phenomenon is localized in the domain interior, periodic boundary conditions effectively emulate an infinite spatial domain under this configuration. When the spatial location of the singularity remains distant from domain boundaries, its interaction with these artificial edges becomes negligible.
\end{remark}

Then at $t_0=0$, we generate $P$ empirical samples $\{ \widetilde{X}^{p}_{0}\}_{p=1:P}$ according to the initial condition of $\widetilde{\rho}_0$ and set up $\widetilde{\alpha}_{0;\mathbf{j}}$ using the Fourier series of $\widetilde{c}_0$. For ease of presenting our algorithm, with a slight abuse of notation, we use $\widetilde{\rho}_n={\frac {M_0}{P}}\, \sum_{p=1}^{P} \delta (x - \widetilde{X}^{p}_{n})$, and 
 \begin{equation}\label{c spectral}
 \widetilde{c}_n=\frac{1}{L^3}\sum_{\mathbf{j} \in \mathcal{H}}\, \widetilde{\alpha}_{n;\mathbf{j}}\, \exp(i 2\pi j_1 \, x_1/L) \exp(i 2\pi j_2\, x_2/L)\exp(i2\pi j_3\, x_3/L)
 \end{equation}
to represent density $\widetilde{\rho}$ and chemical concentration $\widetilde{c}$ at time $t_n$.

Considering the time-stepping for the system (\ref{K-S system}) from $t_n$ to $t_{n+1}$, with $\widetilde{\rho}_n$ and $\widetilde{c}_{n-1}$ known, our algorithm, inspired by the operator splitting technique, consists of two sub-steps: updating chemical concentration $\widetilde{c}$ and updating organism density $\widetilde{\rho}$.

\paragraph{Updating chemical concentration $\widetilde{c}$}
Let $\delta t=t_{n+1}-t_n>0$ be the time step. We discretize the $\widetilde{c}$ equation of \eqref{K-S system} in time by an implicit Euler scheme 
\begin{align}
\epsilon \, ( \widetilde{c}_n - \widetilde{c}_{n-1} )/\delta t
= (\Delta - \lambda^2 )\, \widetilde{c}_n + \widetilde{\rho}_n. \label{Euler4c} 
\end{align}
From Eq.~\eqref{Euler4c}, we obtain the explicit formula for $\widetilde{c}_n$ as follows:
\begin{align}
(\Delta - \lambda^2 -\epsilon/\delta t )\, \widetilde{c}_n
= - \epsilon \, \widetilde{c}_{n-1}/\delta t - \widetilde{\rho}_n.  \label{Euler4c2}
\end{align} 
It follows that 
\begin{align}
	\widetilde{c}_n=\widetilde{c}(\mathbf{x},t_n) &= - \mathcal{K}_{\epsilon,\delta t} \ast (\epsilon \, \widetilde{c}_{n-1}/\delta t + \widetilde{\rho}_n) 
	= - \mathcal{K}_{\epsilon,\delta t} \ast (\epsilon \, \widetilde{c}(\mathbf{x},t_{n-1})/\delta t +\widetilde{\rho}(\mathbf{x},t_n)),
	\label{IPS_Cstep}
\end{align}
where $\mathcal{K}_{\epsilon,\delta t}$ is the Green's function of the operator $\Delta - \lambda^2 -\epsilon/\delta t$ and $\ast$ represents an approximation of spatial convolution that differs from the continuous setup, as $\widetilde{c}$ is computed using truncated Fourier basis functions and $\widetilde{\rho}$ is given by a discrete particle representation.
Unless otherwise stated, all subsequent norms $\|\cdot\|$ will refer to the $L^2$ norms. In the case of $\mathbb{R}^3$, the Green's function $\mathcal{K}_{\epsilon,\delta t}$ reads as follows:
\begin{align}
 \mathcal{K}_{\epsilon,\delta t}
=  \mathcal{K}_{\epsilon,\delta t}(\mathbf{x})= -  
\frac{\exp\{ - \beta \|\mathbf{x} \| \}}{ 4 \pi \|\mathbf{x}\|}, \quad
\beta = \sqrt{\lambda^2 + \epsilon/\delta t}. \label{GreenFunction3D}
\end{align}
Green's function admits a closed-form Fourier transform,
\begin{align}
    \mathcal{F}\mathcal{K}_{\epsilon,\delta t} (\mathbf{\omega})=-\frac{1}{\|\mathbf{\omega}\|^2+\beta^2}.\label{GreenFunction3D_Fourier}
\end{align}
For the term $-\mathcal{K}_{\epsilon,\delta t} \ast  \widetilde{c}_{n-1}$ in Eq.~\eqref{IPS_Cstep}, by Eq.~\eqref{GreenFunction3D_Fourier} it is equivalent to modify Fourier coefficients $\widetilde{\alpha}_{\mathbf{j}}$ to $\widetilde{\alpha}_{\mathbf{j}}/(4\pi^2j_1^2/L^2+4\pi^2j_2^2/L^2+4\pi^2j_3^2/L^2+\beta^2)$.  
 
For the second term $\mathcal{K}_{\epsilon,\delta t} \ast \widetilde{\rho}$, we first approximate $\mathcal{K}_{\epsilon,\delta t}$ using a cosine series expansion. Then, we use the particle representation of $\widetilde{\rho}$ given in Eq.~\eqref{rho_particle} to derive
\begin{align}
    (\mathcal{K}_{\epsilon,\delta t} \ast \widetilde{\rho})_{\mathbf{j}} \approx \frac{M_0}{P}\sum_{p=1}^P \! \frac{\exp(-i2\pi j_1 \widetilde{X}^p_{n;1}/L\!-\!i2\pi j_2 \widetilde{X}^p_{n;2}/L\!-\!i2\pi j_3 \widetilde{X}^p_{n;3}/L)(-1)^{j_1+j_2+j_3}}{4\pi^2j_1^2/L^2+4\pi^2j_2^2/L^2+4\pi^2j_3^2/L^2+\beta^2},
\end{align}
\by{where the factor $(-1)^{j_1+j_2+j_3}$ arises from the redefinition of the Fourier computational domain from $[0, L]^3$ to $[-L/2, L/2]^3$.}

Finally, we summarize the one-step update of the Fourier coefficients of the chemical concentration $\widetilde{c}$ in Algorithm~\ref{alg:SIPF-C}, which follows the same procedure as in the original SIPF method \cite{wang2025novel}. 

\begin{algorithm}[htbp]
    \caption{One-step update of chemical concentration in the SIPF-$r$ method}
    \label{alg:SIPF-C}
    \begin{algorithmic}[1]
        \REQUIRE 
            Distribution $\widetilde{\rho}_n$ represented by empirical samples $\widetilde{X}_n$,  
            initial concentration $\widetilde{c}_{n-1}$ represented by Fourier coefficients $\widetilde{\alpha}_{n-1}$.
        \FOR{$\mathbf{j} \in \mathcal{H}$}
            \STATE $\widetilde{\alpha}_{n;\mathbf{j}} \gets \dfrac{\epsilon \widetilde{\alpha}_{n-1;\mathbf{j}}}{\delta t(4\pi^2 j_1^2/L^2 + 4\pi^2 j_2^2/L^2 + 4\pi^2 j_3^2/L^2 + \beta^2)}$
            \STATE $F_{\mathbf{j}} \gets 0$
            \FOR{$p = 1$ \TO $P$}
                \STATE $F_{\mathbf{j}} \gets F_{\mathbf{j}} + \exp(-i2\pi j_1 \widetilde{X}^p_{n;1}/L-i2\pi j_2 \widetilde{X}^p_{n;2}/L-i2\pi j_3 \widetilde{X}^p_{n;3}/L)$
            \ENDFOR
            \STATE $F_{\mathbf{j}} \gets F_{\mathbf{j}} \cdot \dfrac{(-1)^{j_1+j_2+j_3}}{4\pi^2j_1^2/L^2+4\pi^2j_2^2/L^2+4\pi^2j_3^2/L^2+\beta^2} \cdot \dfrac{M_0}{P}$
        \ENDFOR
        \STATE $\widetilde{\alpha}_n \gets \widetilde{\alpha}_n + F$
        \ENSURE
            Updated chemical concentration field from $\widetilde{c}_{n-1}$ to $\widetilde{c}_n$ via $\widetilde{\alpha}_n$.
    \end{algorithmic}
\end{algorithm}

\paragraph{Updating density of active particles $\widetilde{\rho}$}
 In the one-step update of density $\widetilde{\rho}_n$ represented by particles $\{\widetilde{X}_n^p\}_{p=1:P}$, we apply the Euler-Maruyama scheme to solve the stochastic differential equation (SDE) 
\begin{align}
	\widetilde{X}^{p}_{n+1} = \widetilde{X}^{p}_{n} 
	+\chi \nabla_\mathbf{x} \widetilde{c}(\widetilde{X}^{p}_{n},t_n)\delta t
	+ \sqrt{2\,\mu\, \delta t} \, N^p_n, 
	 \label{OnestepX_EM}
\end{align}
where the variables $N^p_n$ are i.i.d. standard normal random variables corresponding to the Brownian paths in the SDE formulation. 
For $n>1$, substituting Eq.~\eqref{IPS_Cstep} in Eq.~\eqref{OnestepX_EM} gives:
\begin{align}
	\widetilde{X}^{p}_{n+1} = \widetilde{X}^{p}_{n} 
	-\chi \nabla_\mathbf{x} \mathcal{K}_{\epsilon,\delta t} \ast (\epsilon\, \widetilde{c}_{n-1}(\mathbf{x})/\delta t + \widetilde{\rho}_n(\mathbf{x}))|_{\mathbf{x}=\widetilde{X}^{p}_{n}}\delta t
	+ \sqrt{2\,\mu\, \delta t} \,  N^p_n, 
	\label{IPS_Xstep_full}
\end{align}
from which $\widetilde{\rho}_{n+1}(\mathbf{x})$ is constructed via Eq.~\eqref{rho_particle}. 

In this particle formulation, the computation of the spatial convolution differs slightly from that in the update of $\widetilde{c}$ (i.e., Eq.~\eqref{IPS_Cstep}).

For the term $\nabla_{\mathbf{x}}\mathcal{K}_{\epsilon,\delta t} \ast \, \widetilde{c}_{n-1}(\widetilde{X}^p_n)$, to avoid the singular points of $\nabla_{\mathbf{x}}\mathcal{K}_{\epsilon,\delta t}$, we evaluate the integral with quadrature points that are away from $0$. Precisely, denote the standard quadrature point in $\Omega$ as 
\begin{equation}
x_{\mathbf{j}}=( j_1\, L/H,j_2\, L/H, j_3\, L/H), \label{grid}
\end{equation}
where $j_1$, $j_2$, $j_3$ are integers ranging from $-H/2$ to $H/2-1$. When computing $\nabla_{\mathbf{x}}\mathcal{K}_{\epsilon,\delta t} \ast \, \widetilde{c}_{n-1}(\widetilde{X}^p_n)$, we evaluate $\nabla_{\mathbf{x}}\mathcal{K}_{\epsilon,\delta t}$ at $\{\widetilde{X}^p_n+\bar{X}^p_n-x_{\mathbf{j}}\}_{\mathbf{j}}$, where a small spatial shift is defined as $\bar{X}^p_n=\frac{L}{2H}+\lfloor \frac{\widetilde{X}^p_n}{L/H}\rfloor \frac{L}{H}-\widetilde{X}_n^p$ and $\widetilde{c}$ at $\{x_{\mathbf{j}}-\bar{X}^p_n\}_{\mathbf{j}}$ correspondingly. The latter is computed by inverse Fourier transform of the shifted coefficients, with $\widetilde{\alpha}_{\mathbf{j}}$ modified to $\widetilde{\alpha}_{\mathbf{j}}\exp(-i2\pi j_1\bar{X}^p_{n;1}/L- i2\pi j_2\bar{X}^p_{n;2}/L-i2\pi j_3\bar{X}^p_{n;3}/L) $, where $(\bar{X}^p_{n;i})$ denotes the $i$-th component of $\bar{X}^p_{n}$.

Motivated by mini-batch sampling \cite{DL_2016,sommerfeld2019optimal,DP_22,wang2024deepparticle,xie2024randomized,bidir-DP-2026} and the random batch method (RBM) \cite{jin2020random, jin2021random, cai2024convergence, huang2025mean}, for each particle $\widetilde{X}_n^p$, we choose a small batch $C_p$ of size $R$ randomly with replacement. We restrict the interaction of $\widetilde{X}_n^p$ to particles within this batch; specifically, we approximate the interaction term $\chi \delta t \nabla_{\mathbf{x}}\mathcal{K}_{\epsilon,\delta t} \ast \widetilde{\rho}(\widetilde{X}^p_n,t_{n})$ by the random batch sum $\sum_{s \in C_p, s \neq p} \frac{\chi M_0 \delta t}{R} \nabla_{\mathbf{x} } \mathcal{K}_{\epsilon, \delta t}(\widetilde{X}^p_n - \widetilde{X}^s_n)$.

We summarize the one-step update (for $n>1$) of the density in the SIPF-$r$ method as in Algorithm~\ref{alg:SIPF-$r$ho}.

\begin{algorithm}[htbp]
    \caption{One-step update of density in the SIPF-$r$ method}
    \label{alg:SIPF-$r$ho}
    \begin{algorithmic}[1]
        \REQUIRE 
            Distribution $\widetilde{\rho}_n$ represented by empirical samples $\widetilde{X}_n$, concentration $\widetilde{c}_{n-1}$ represented by Fourier coefficients $\widetilde{\alpha}_{n-1}$.
        \FOR{$p=1$ \TO $P$}
            \STATE $\widetilde{X}_{n+1}^p \gets \widetilde{X}_{n}^p + \sqrt{2\mu\delta t}N^p_n$ \COMMENT{$N^p_n$ is a standard normal random variable}
            \STATE $C_p \gets$ random subset of $\{1, \dots, P\}$ with replacement, size $R$
            \STATE $\widetilde{X}_{n+1}^p \gets \widetilde{X}_{n+1}^p - \sum_{s \in C_p, s \neq p} \frac{\chi M_0 \delta t}{R} \nabla_{\mathbf{x}} \mathcal{K}_{\epsilon, \delta t}(\widetilde{X}^p_n - \widetilde{X}^s_n)$
            \STATE $\bar{X}^p_n \gets \frac{L}{2H} + \lfloor \frac{\widetilde{X}^p_n}{L/H} \rfloor \frac{L}{H} - \widetilde{X}_n^p$
            \FOR{$(\mathbf{j}) \in \mathcal{H}$}
                \STATE $F_{\mathbf{j}} \gets \nabla_{\mathbf{x}} \mathcal{K}_{\epsilon,\delta t}(\widetilde{X}^p_n+\bar{X}_n^p-x_{\mathbf{j}})$ \quad \COMMENT{$x_{\mathbf{j}}$ from Eq.~\eqref{grid}}
                \STATE $G_{\mathbf{j}} \gets \widetilde{\alpha}_{n-1; \mathbf{j}} \exp(-i2\pi j_1\bar{X}^p_{n;1}/L- i2\pi j_2\bar{X}^p_{n;2}/L-i2\pi j_3\bar{X}^p_{n;3}/L)$
            \ENDFOR
            \STATE $\check{G} \gets \text{iFFT}(G)$
            \STATE $\widetilde{X}_{n+1}^p \gets \widetilde{X}_{n+1}^p - \epsilon \chi(F,\check{G}) \frac{L^3}{H^3}$ \quad \COMMENT{$(\cdot,\cdot)\frac{L^3}{H^3}$ denotes an inner product corresponding to $L^2(\Omega)$ quadrature}
        \ENDFOR
        \ENSURE
            Updated distribution $\widetilde{\rho}_{n+1}$ represented by $\widetilde{X}_{n+1}$.
    \end{algorithmic}
\end{algorithm}

\by{\begin{remark}
The spatial shift $\bar{X}_n^p$ acts as a mild perturbation to avoid the singular attractive force induced by the Green's function $\mathcal{K}_{\epsilon,\delta t}$. While it may technically introduce bias, due to the even structure of $\mathcal{K}_{\epsilon,\delta t}$ and the regularity of $\widetilde{c}$, the influence of this shift is small when evaluating $(F,\check{G})$. Consequently, particle aggregation remains driven by the intrinsic nonlinear dynamics rather than the grid structure. Combined with the stochasticity of the random batch method, this ensures that no systematic directional bias is introduced, which allows for robust singularity detection independent of specific grid points.
\end{remark}}

Combining Eq.~\eqref{IPS_Cstep} and Eq.~\eqref{IPS_Xstep_full}, we conclude that the recursion from \\
$(\{\widetilde{X}^{p}_{n}\}_{p=1:P}, \widetilde{\rho}_n(\mathbf{x}), \widetilde{c}_{n-1}(\mathbf{x}))$  
to $(\{\widetilde{X}^{p}_{n+1}\}_{p=1:P}, \widetilde{\rho}_{n+1}(\mathbf{x}),\widetilde{c}_n(\mathbf{x}))$  
is thus fully defined. We summarize the SIPF-$r$ method in the following Algorithm~\ref{alg:SIPF}.

\begin{algorithm}[htbp]
    \caption{Stochastic Interacting Particle-Field Method}
    \label{alg:SIPF}
    \begin{algorithmic}[1]
        \REQUIRE 
            Initial distribution $\rho_0$, initial concentration $c_0$.
        \STATE Generate $P$ i.i.d.\ samples according to distribution $\rho_0$: $\widetilde{X}_0^1, \widetilde{X}_0^2, \dots, \widetilde{X}_0^P$.
        \FOR{$p = 1$ \TO $P$}
            \STATE Compute $\widetilde{X}^{p}_{1}$ by Eq.~\eqref{OnestepX_EM}, with $c_{-1} = c_0$.
        \ENDFOR
        \STATE Compute $\widetilde{c}_1$ by Alg.~\ref{alg:SIPF-C} with $c_0$ and $\widetilde{\rho}_1 = \sum_{p=1}^P \frac{M_0}{P} \delta_{\widetilde{X}^p_1}$.
        \FOR{$n = 2$ \TO $N = \lfloor T / \delta t \rfloor$}
            \STATE Compute $\widetilde{X}_{n}$ by Alg.~\ref{alg:SIPF-$r$ho} with $\widetilde{\rho}_{n-1}$ and $\widetilde{c}_{n-2}$.
            \STATE Compute $\widetilde{c}_{n}$ by Alg.~\ref{alg:SIPF-C} with $\widetilde{c}_{n-1}$ and $\widetilde{\rho}_{n} = \sum_{p=1}^P \frac{M_0}{P} \delta_{\widetilde{X}^p_n}$.
        \ENDFOR
        \ENSURE
            Final particle distribution $\widetilde{\rho}_N$ and concentration field $\widetilde{c}_N$.
    \end{algorithmic}
\end{algorithm}

\paragraph{Computational Complexity}
We briefly analyze the complexity of the proposed SIPF-$r$ method. The memory usage is $\mathcal{O}(P + |\mathcal{H}|)$, where $|\mathcal{H}|$ denotes the total number of Fourier modes. Regarding the computational cost, the original SIPF method typically requires $\mathcal{O}(P^2)$ operations for pairwise particle interactions. In contrast, by incorporating the RBM, the interaction cost in the SIPF-$r$ method is reduced to $\mathcal{O}(PR)$, where $R$ is the batch size.

The overall per-step complexity of Algorithm~\ref{alg:SIPF} is $\mathcal{O}(PR + P |\mathcal{H}|\log(|\mathcal{H}|))$. Crucially, the error analysis in Theorem \ref{theorem} reveals that the error introduced by the random batch approximation scales only as $\mathcal{O}(\sqrt{\delta t/R})$. This allows us to choose $R \ll P$ (e.g., $R=100$ vs. $P=10^4$) to achieve a substantial speedup while maintaining the desired accuracy.

\paragraph{Particle-wise Independence due to RBM}

In the above derivation, $\{\widetilde{X}^{p}_{n}\}_{p=1:P}$ are i.i.d.\ samples with distribution $\widetilde{\rho}_n$ and independent of $\widetilde{c}_{n-1}$. The one-step trajectories follow the discrete-time rule:
\begin{equation}\label{X SIPF}
    \widetilde{X}_{t_{n+1}} = \widetilde{X}_{t_n} + \chi \nabla \widetilde{c}(\widetilde{X}_{t_n}, t_n) \delta t + \int_{t_n}^{t_{n+1}} \sqrt{2\mu} \, dW_s,
\end{equation}
where $\nabla \widetilde{c}$ is computed via Eq.~\eqref{IPS_Cstep}, and $W_s$ denotes the Brownian motion. It is worth noting that, for the updated position of the $p$-th particle $\widetilde{X}_{n+1}^p$ by Eq.~\eqref{OnestepX_EM}, the interaction term, $\nabla_{\mathbf{x}}\mathcal{K}_{\epsilon,\delta t} \ast \, \widetilde{\rho}(\widetilde{X}^p_n,t_{n})$ is computed by $\sum_{s \in C_p, s \neq p} \frac{\chi M_0 \delta t}{R} \nabla_{\mathbf{x} } \mathcal{K}_{\epsilon, \delta t}(\widetilde{X}^p_n - \widetilde{X}^s_n)$, where the selection of $C_p$ is independent of $\widetilde{X}^p_n$ and hence $\{\widetilde{X}^s_n\}_{s\in C_p}$ can be viewed as i.i.d. samples of $\widetilde{\rho}_n$ independent of $\widetilde{c}_{n-1}$ and $\widetilde{X}^p_n$. Together with the independent Brownian motion term, we can deduce the independence of $\{\widetilde{X}_{n+1}^p\}_{p=1:P}$.

Correspondingly, we denote the exact dynamics of the system by $X_{t}$, a $\rho(\cdot, t)$-distributed random variable evolving continuously in time:
\begin{equation}\label{X exact}
X_{t} = X_{t_0} + \chi \int_{t_0}^{t} \nabla c(X_s, s) \, ds + \int_{t_0}^{t} \sqrt{2\mu} \, dW_s, \quad {X}_{t_0}=\widetilde{X}_{t_0},
\end{equation}
where $c(\cdot, s)$ is the exact concentration field, and the integral describes how the gradient field evolves in continuous time. Both processes share the same Brownian motion $W_s$, indicating that both processes are driven by the same source of randomness.

\section{$L^2$ Convergence of the SIPF-$r$ Method to Smooth Solutions}\label{section: proof}

We now prove the convergence of the SIPF-$r$ method to classical solutions of the 3D parabolic-parabolic Keller-Segel equations. To ensure the validity of the following analysis, we introduce a set of assumptions that impose structure on the concentration fields and their gradients.

\subsection{Assumptions and Theorem Statement}
To establish the convergence analysis of the SIPF-$r$ method, we first introduce assumptions that ensure the boundedness of approximation errors for particles and gradients of the chemical concentration \by{at any time within the temporal computational domain $[0, T]$}. Specifically, we make the following assumptions.

\begin{assumption}\label{assump:1}
There exist constants $ M_1, M_2 > 0 $ such that for all $t \in [0, T]$ and $\mathbf{x} \in \mathbb{R}^3$,
\begin{align}
    \|\widetilde{X}_{t} - X_{t}\| & \leq M_1, \label{X～-X} \\
    \|\nabla \widetilde{c}(\mathbf{x}, t)-\nabla c(\mathbf{x}, t)\| & \leq M_2. \label{gradient c~-c}
\end{align}
\end{assumption}
We note that this assumption only requires the errors to be bounded, but does not require them to converge to zero. The convergence of these errors to zero will be established in the subsequent analysis.
Furthermore, the boundedness condition in Eq.~\eqref{X～-X} can be derived from Eqs.~\eqref{X SIPF}-\eqref{X exact} and  Assumption~\ref{assump:2}(c); Eq.~\eqref{gradient c~-c} follows directly from the uniform bound in Assumption \ref{assump:2}(c) detailed later. 

\begin{assumption}\label{assump:2}
Suppose both $\nabla \widetilde{c}$ and $\nabla c$ satisfy Lipschitz continuity conditions in space and time, along with regularity and boundedness properties as follows:

\noindent\textbf{(a) (Spatial Lipschitz Continuity)} There exists a constant $K > 0$, depending on the regularity of $\nabla \widetilde{c}$ and $\nabla c$, as well as the parameters $\epsilon$ and $\lambda$ in the system \eqref{K-S system}, such that for all $t \in [0, T]$ and $\mathbf{x}, \mathbf{y} \in \mathbb{R}^3$,
\begin{equation}
\max\big(\|\nabla \widetilde{c}(\mathbf{x}, t) - \nabla \widetilde{c}(\mathbf{y}, t)\|,\,
            \|\nabla c(\mathbf{x}, t) - \nabla c(\mathbf{y}, t)\| \big)
\leq K \|\mathbf{x} - \mathbf{y}\|.
\end{equation}
    This implies that the second derivatives (Hessian entries) $\nabla^2 \widetilde{c}(\mathbf{x}, t)$ exist almost everywhere and satisfy:
    \begin{equation}
    \sup\limits_{\mathbf{x}\in \mathbb{R}^3, t \in [0, T]}\big(\|\nabla^2 \widetilde{c}(\mathbf{x}, t)\|\big) \leq K.
    \end{equation}

\noindent\textbf{(b) (Temporal Lipschitz Continuity)} There exists a constant $K_1 > 0$, depending on the regularity of $\nabla c$ and the parameters $\epsilon$ and $\lambda$ in the system \eqref{K-S system}, such that for any $t_1, t_2 \in [0, T]$ and $\mathbf{x} \in \mathbb{R}^3$,
    \begin{equation}
    \|\nabla c(\mathbf{x}, t_1) - \nabla c(\mathbf{x}, t_2)\| \leq K_1 |t_1 - t_2|.
    \end{equation}

\noindent\textbf{(c) (Uniform Boundedness)} There exists a constant $M_3 > 0$, depending on the regularity of $\nabla c$ and the parameters $\epsilon$ and $\lambda$ in the system \eqref{K-S system}, such that for all $t \in [0, T]$ and $\mathbf{x} \in \mathbb{R}^3$:
    \begin{equation}
    \max\big(\|\nabla c(\mathbf{x}, t)\|, \|\nabla \widetilde{c}(\mathbf{x}, t)\|\big) \leq M_3.
    \end{equation}

\noindent\textbf{(d) (Regularity of Time Derivatives)} The exact solution $c(\mathbf{x}, t)$ is assumed to be sufficiently smooth in time and space such that both $\partial_{t}^2 c(\mathbf{x}, t)$ and $\nabla \partial_{t}^2 c(\mathbf{x}, t)$ are bounded. There exists a constant $K_2 > 0$ such that for all $t \in [0, T]$:
    \begin{equation}
    \max\left(\sup_{\mathbf{x} \in \mathbb{R}^3} \|\partial_{t}^2 c(\mathbf{x}, t)\|, \sup_{\mathbf{x} \in \mathbb{R}^3} \|\nabla \partial_{t}^2 c(\mathbf{x}, t)\|\right) \leq K_2.
    \end{equation}
\end{assumption}

\begin{remark}
\by{These regularity and boundedness conditions reflect inherent properties of the SIPF-$r$ method rather than external constraints on the KS system. As analyzed in Subsection~\ref{subsec:stability}, the discrete Fourier representation and Brownian diffusion naturally impose spectral regularity and 
suppress singular clustering, ensuring $\nabla\widetilde{c}$ remains well-behaved throughout the computation. The boundedness conditions in Assumption \ref{assump:2} are validated through numerical simulations in Subsection \ref{subsection: Accuracy of SIPF-$r$ Method}.}
\end{remark}

\by{Assumptions~\ref{assump:1} and \ref{assump:2} characterize a regular regime where the exact solution remains smooth over the time interval $[0,T]$.} We now state our main theorem, which quantifies the convergence of the SIPF-$r$ method on $[0,T]$.

\by{\begin{theorem}\label{theorem}
Suppose that the exact solutions $(\rho,c)$ and the numerical solutions $(\widetilde{\rho},\widetilde{c})$ obtained by the SIPF-$r$ method in $\mathbb{R}^3$ satisfy Assumptions \ref{assump:1} and \ref{assump:2} uniformly for all $t\in[0,T]$. Let $ H $, $ P $, $ R $, and $ \delta t $ denote the number of Fourier modes, the number of particles, the batch size, and the uniform time step in the SIPF-$r$ method, respectively. Then, for all discrete time levels $ t_n = n\delta t \leq T $, the following error estimates hold with high probability:

For all $ n \in \{0, 1, \dots, \lfloor T/\delta t\rfloor\}$, the $1$-Wasserstein distance (defined in Eq.~\eqref{def:wasserstein distance}) between $ \widetilde{\rho}_{t_n}$ and $ \rho_{t_n} $ satisfies $\mathcal{W}_1(\widetilde{\rho}_{t_n}, \rho_{t_n})=\mathcal{O}\left(\frac{1}{H^2} + \sqrt{\frac{1}{P}}  + \sqrt{\delta t}\right)$, and the maximum error in the truncated Fourier coefficients of $\widetilde{c}_{t_n}$ and $c_{t_n}$ satisfies that $\max_{\mathbf{j} \in \mathcal{H}}|\widetilde{\alpha}_{t_n;\mathbf{j}} - \alpha_{t_n;\mathbf{j}}|=\mathcal{O}\left(\frac{1}{H} + \frac{H}{\sqrt{P}}  + H\delta t\right)$. More specifically, for all $ n \in \{0, 1, \dots, \lfloor T/\delta t\rfloor\} $, the errors are bounded with high probability by:
\begin{equation}\label{theorem inequality}
\begin{aligned}
\mathcal{W}_1(\widetilde{\rho}_{t_n}, \rho_{t_n}) &\leq \frac{S_0}{H^2} + S_1\sqrt{\delta t} + \frac{S_2}{\sqrt P} + \frac{S_3\sqrt{\delta t}}{\sqrt{R}}, \\
\max_{\mathbf{j} \in \mathcal{H}}|\widetilde{\alpha}_{t_n;\mathbf{j}} - \alpha_{t_n;\mathbf{j}}| &\leq \frac{S_4}{H} + S_5H\delta t + \frac{S_6H}{\sqrt{P}} + \frac{S_7H\delta t}{\sqrt{R}},
\end{aligned}
\end{equation}
where $ S_i $, $ i = 0, \dots, 7$, denote positive constants specified in Eqs.~\eqref{eq:wasserstein_bound}-\eqref{general alpha}. 
\end{theorem}}

A direct consequence of Theorem \ref{theorem} is that the SIPF-$r$ solution converges to the exact solution as $\delta t \to 0$ and $H, P \to \infty$. 
Specifically, the density $\widetilde{\rho}_{t_n}$ converges in the 1-Wasserstein metric even without imposing specific scaling relations between the discretization parameters, as all error terms in the first inequality of Eq.~\eqref{theorem inequality} naturally vanish. The convergence of the chemical concentration $\widetilde{c}_{t_n}$ requires the scaling conditions $H/\sqrt{P} \to 0$ and $H\delta t \to 0$ to control the statistical and discretization errors, respectively.

\begin{remark}[Practical convergence and scaling]\label{remark:practical}
\by{Theorem~\ref{theorem} establishes error estimates up to time $T$, the maximal time at which Assumptions~\ref{assump:1}-\ref{assump:2} hold. Beyond $T$, although the theoretical bounds no longer apply, the Fourier spectral cutoff $H$ and the parabolic structure of the equation for $c$ act as a natural low-pass filter that prevents the numerical solution $\widetilde{c}$ from blowing up. Note that this numerical upper bound grows with $H$. By repeating the experiment with increasing values of $H$, we obtain a reliable detection of the time and location of the singularity.  

Regarding parameter scaling, the stochastic error term $\mathcal{O}(\sqrt{\delta t/R})$ is negligible relative to the leading-order discretization errors. Hence, a moderate batch size (e.g., $R=100$) is sufficient in practice, as the total error is dominated by $H$ and $\delta t$. Furthermore, while the theoretical gradient error bound scales with $H$ due to differentiation, the particle density error benefits from the smoothing effect of time integration. Numerical experiments in Section~\ref{section: experiment} verify that moderate parameters ($H=24$, $P=10^4$, and $\delta t=10^{-4}$) are sufficient to guarantee convergence and blow-up detection for our numerical examples, which demonstrates the efficiency of the proposed method.}
\end{remark}

\subsection{Several lemmas}
The result of Theorem \ref{theorem} relies on the following lemmas concerning the change in single-step update error of the SIPF-$r$ method. We first prove several lemmas in this subsection and leave the complete proof of Theorem \ref{theorem} to the next subsection. 


From Eqs.~\eqref{trig_ser}-\eqref{c fourier}, the error between $\widetilde{c}$ and $c$ can be decomposed into two components: the error in their Fourier coefficients and the truncation error of $c$. As the number of Fourier modes $H$ tends to infinity, and given the smoothness of $c$, the truncation error becomes negligible and can be omitted from the analysis. We now focus on the error analysis between the Fourier coefficients $\widetilde{\alpha}_{\mathbf{j}}$ and $\alpha_{\mathbf{j}}$ of $\widetilde{c}$ and $c$, as presented in the following lemma.

\begin{lemma}\label{lemma1}
For all $ n \in \mathbb{N_{+}}$ and all $ \mathbf{j} \in \mathcal{H}$ (the same index set as in Eq.~\eqref{index set H}), under Assumptions \ref{assump:1} and \ref{assump:2}, the following inequality holds with high probability:
\begin{align}
|\widetilde{\alpha}_{t_n;\mathbf{j}} - \alpha_{t_n;\mathbf{j}}| 
\leq& |\widetilde{\alpha}_{t_{n-1};\mathbf{j}} - \alpha_{t_{n-1};\mathbf{j}}| + C_1\delta t^2 + \frac{C_2\|\omega_{\mathbf{j}}\|\delta t}{\sqrt{P}}\notag + C_3\|\omega_{\mathbf{j}}\| \delta t \mathbb{E}[\|\widetilde{X}_{t_n} - X_{t_n}\|],
\end{align}
where $ C_1, C_2, C_3$ are constants, and $t_n = n \delta t$. 
\end{lemma}

\begin{proof} 
We first write the frequency $\omega_{\mathbf{j}} = \left(\frac{2\pi j_1}{L}, \frac{2\pi j_2}{L}, \frac{2\pi j_3}{L}\right)$, and define the notation  
$Z_{\mathbf{j}} = \left(\|\omega_{\mathbf{j}}\|^2 + \lambda^2\right) \cdot \frac{\delta t}{\epsilon}$, which satisfies $Z_{\mathbf{j}} \to 0$ as $\delta t \to 0$. Recall the update formula for the numerical approximation of the Fourier coefficients $\widetilde{\alpha}_{t_n;\mathbf{j}}$ (from Eq.~\eqref{IPS_Cstep}):
\begin{equation}\label{num_recurrence}
(1 + Z_{\mathbf{j}}) \widetilde{\alpha}_{t_n;\mathbf{j}} = \widetilde{\alpha}_{t_{n-1};\mathbf{j}} + \frac{\delta t}{\epsilon}\mathcal{F}_{\mathbf{j}}[\widetilde{\rho}(\mathbf{x}, t_n)],
\end{equation}
where $\mathcal{F}_\mathbf{j}[\widetilde{\rho}(\mathbf{x}, t_n)]=\frac{M_0}{P} \sum_{p=1}^P e^{-i\mathbf{\omega}_{\mathbf{j}} \cdot \widetilde{X}^p_{t_n}}$ represents the Fourier coefficient of $\widetilde{\rho}(\mathbf{x}, t_n)$ at the frequency $\omega_{\mathbf{j}}$.

The exact solution $c(\mathbf{x}, t)$ satisfies the continuous equation:
\begin{equation}
\epsilon \frac{\partial c}{\partial t} = \Delta c - \lambda^2 c + \rho.
\end{equation}
Integrating this equation from $t_{n-1}$ to $t_n$ and applying the Taylor expansion with integral remainder to the time derivative term yields $c(t_{n-1}) = c(t_n) - \delta t \partial_t c(t_n) + \int_{t_{n-1}}^{t_n} (s - t_{n-1}) \partial_{t}^2 c(s) \, ds$.  This allows us to express the exact solution in a form compatible with the implicit Euler scheme:
\begin{equation}
\frac{c(t_n) - c(t_{n-1})}{\delta t} = \frac{1}{\epsilon}(\Delta c(t_n) - \lambda^2 c(t_n) + \rho(t_n)) + R_n(\mathbf{x}).
\end{equation}

By Assumption \ref{assump:2}(d), $\|\partial_{t}^2 c\|$ is bounded by $K_2$. The temporal truncation error $R_n(\mathbf{x})$ satisfies
\begin{equation}\label{R_n}
    \|R_n\| \leq \frac{\delta t}{2} \, \sup_{s \in [t_{n-1}, t_n]} \|\partial_{t}^2 c(\cdot, s)\| \leq C_1 \, \delta t,
\end{equation}
where $C_1=\frac{K_2}{2}$ is a constant. Taking the Fourier transform of the discrete relation for the exact solution and rearranging terms, we get
\begin{equation}\label{exact_recurrence}
(1 + Z_{\mathbf{j}}) \alpha_{t_n;\mathbf{j}} = \alpha_{t_{n-1};\mathbf{j}} + \frac{\delta t}{\epsilon}\mathcal{F}_{\mathbf{j}}[\rho(t_n)] + \delta t \mathcal{F}_{\mathbf{j}}[R_n].
\end{equation}

Combining Eq.~\eqref{exact_recurrence} and Eq.~\eqref{num_recurrence}, we obtain
\begin{align}\label{tilde_alpha-alpha}
|\alpha_{t_n;\mathbf{j}} \!-\! \widetilde{\alpha}_{t_n;\mathbf{j}}| =& \frac{1}{1+Z_{\mathbf{j}}} \left|\left(\alpha_{t_{n-1};\mathbf{j}} - \widetilde{\alpha}_{t_{n-1};\mathbf{j}}\right) \!+\! \frac{\delta t}{\epsilon} \left(\mathcal{F}_{\mathbf{j}}[\rho(t_n)] - \mathcal{F}_{\mathbf{j}}[\widetilde{\rho}_n]\right) \!+\! \delta t\mathcal{F}_{\mathbf{j}}[R_n]\right|\notag \\
\leq& |\alpha_{t_{n-1};\mathbf{j}} - \widetilde{\alpha}_{t_{n-1};\mathbf{j}}|+\frac{\delta t}{\epsilon}|\mathcal{F}_{\mathbf{j}}[\rho(t_n)] - \mathcal{F}_{\mathbf{j}}[\widetilde{\rho}_n]|+\delta t|\mathcal{F}_{\mathbf{j}}[R_n]|.
\end{align}

To bound the term $|\mathcal{F}_{\mathbf{j}}[\rho(t_n)] - \mathcal{F}_{\mathbf{j}}[\widetilde{\rho}_n]|$, we first state a generalization of the mean value theorem to complex-valued functions.

Let $ G $ be an open subset of $ \mathbb{R}^3$, and let $ f: G \to \mathbb{C} $ be a continuously differentiable function on $G$. Fix points $ \mathbf{x}, \mathbf{y} \in G $ such that the line segment connecting $ \mathbf{x} $ and $ \mathbf{y} $ lies entirely within $ G $. There exists $c_1, c_2 \in (0, 1)$ such that
\begin{equation}\label{complex mean value theorem}
f(\mathbf{y}) - f(\mathbf{x}) = \operatorname{Re}\Big(\nabla f((1-c_1) \mathbf{x} + c_1 \mathbf{y})(\mathbf{y} - \mathbf{x})\Big) + i \operatorname{Im}\Big(\nabla f((1-c_2) \mathbf{x} + c_2 \mathbf{y})(\mathbf{y} - \mathbf{x})\Big).
\end{equation}

The proof of \eqref{complex mean value theorem} is direct. First, we define the function
\begin{equation}
g(t) = f((1-t) \mathbf{x} + t \mathbf{y}), \quad t \in [0, 1].
\end{equation}
Since $g$ is also a continuously differentiable function, the mean value theorem implies that there exist points $c_1, c_2 \in (0, 1)$ such that
\begin{equation}
\operatorname{Re}(g'(c_1)) = \operatorname{Re}(g(1) - g(0)),\quad \operatorname{Im}(g'(c_2)) = \operatorname{Im}(g(1) - g(0)),
\end{equation}
which implies Eq.~\eqref{complex mean value theorem}. Applying this result to $f(\mathbf{x})=e^{-i\mathbf{\omega}_{\mathbf{j}}\mathbf{x}}$, we obtain:

\begin{align}\label{exponential X}
|e^{-i\mathbf{\omega}_{\mathbf{j}}\cdot \widetilde{X}^p_{t_n}}-e^{-i\mathbf{\omega}_{\mathbf{j}}\cdot X^p_{t_n}}| &\leq \|\omega_{\mathbf{j}}\cdot \sin(\omega_{\mathbf{j}}((1-c_1)\widetilde{X}^p_{t_n}+c_1X^p_{t_n})) \notag\\
&+i\omega_{\mathbf{j}}\cdot \cos(\omega_{\mathbf{j}}((1-c_2)\widetilde{X}^p_{t_n}+c_2X^p_{t_n}))\|\cdot \|\widetilde{X}^p_{t_n}-X^p_{t_n}\| \notag\\
&\leq \sqrt{2}\|\omega_{\mathbf{j}}\|\cdot \|\widetilde{X}^p_{t_n}-X^p_{t_n}\|.
\end{align}

Based on Eq.~\eqref{exponential X}, we have
\begin{align}
|\mathcal{F}_{\mathbf{j}}[\rho(t_n)] \!-\! \mathcal{F}_{\mathbf{j}}[\widetilde{\rho}_n]|\leq \left\| \frac{M_0}{P}\!\sum_{p=1}^P(e^{-i\mathbf{\omega}_{\mathbf{j}} X^p_{t_n}}\!-\!e^{-i\mathbf{\omega}_{\mathbf{j}} \widetilde{X}^p_{t_n}})\right\| \!\leq\! \sqrt{2}M_0\|\omega_{\mathbf{j}}\|\!\sum_{p=1}^P\!\frac{\|\widetilde{X}^p_{t_n}\!-\!X^p_{t_n}\|}{P}.
\end{align}
Let $Y_p = \|\widetilde{X}_{t_n}^p - X_{t_n}^p\|$, where $\{Y_p\}_{p=1}^P$ are i.i.d. random variables. This follows from the fact that the particles $\{X_{t_n}^p\}_{p=1}^P$ and $\{\widetilde{X}_{t_n}^p\}_{p=1}^P$ are separately i.i.d. Specifically, the i.i.d. property of $\{\widetilde{X}_{t_n}^p\}_{p=1}^P$ is ensured by the RBM described in Alg.~\ref{alg:SIPF-$r$ho}. Based on Assumption \ref{assump:1}, $Y_p$ is bounded. The empirical mean is defined as $\bar{Y}_P = \frac{1}{P} \sum_{p=1}^P Y_p$ and the expectation of $Y_p$ is $\mu = \mathbb{E}[Y_p] = \mathbb{E}[\|\widetilde{X}_{t_n} - X_{t_n}\|]$.
 
According to Bernstein's inequality, for i.i.d. random variables $Y_1, Y_2, \dots, Y_P$ with  $|Y_p - \mu| \leq M_1$ (from Assumption \ref{assump:1}) almost surely, the probability that the empirical mean deviates from the expectation is bounded as
\begin{equation}
\mathbb{P}\left( |\bar{Y}_P - \mu| \geq \eta \right) \leq 2\exp\left(-\frac{P\eta^2}{2\sigma^2 + \frac{2M_1\eta}{3}}\right),
\end{equation}
where $\sigma^2 = \mathbb{E}[(Y_p-\mu)^2]\leq M_1^2$. With high probability (e.g., $1 - \delta_0$ for very small $\delta_0 > 0$), the following estimate holds
\begin{equation}
|\bar{Y}_P - \mu| \leq \sqrt{\frac{2\sigma^2 \ln(2/\delta_0)}{P}} + \frac{2M_1 \ln(2/\delta_0)}{3P}.
\end{equation}
This implies that, with probability $1 - \delta_0$,  
\begin{equation}\label{rho_fourier}
|\mathcal{F}_{\mathbf{j}}[\rho(t_n)] \!-\! \mathcal{F}_{\mathbf{j}}[\widetilde{\rho}_n]| \leq \sqrt{2}M_0\|\omega_{\mathbf{j}}\| \left(\!\mathbb{E}[\|\widetilde{X}_{t_n} - X_{t_n}\|]\!+\!\sqrt{\frac{2\sigma^2 \ln(2/\delta_0)}{P}} \!+\! \frac{2M_1 \ln(2/\delta_0)}{3P}\!\right).
\end{equation}

Combining Eqs.~\eqref{R_n}\eqref{tilde_alpha-alpha}\eqref{rho_fourier} above, we conclude that, with high probability
\begin{align}\label{alpha~-alpha}
|\widetilde{\alpha}_{t_n;\mathbf{j}} - \alpha_{t_n;\mathbf{j}}|
\leq& |\widetilde{\alpha}_{t_{n-1};\mathbf{j}} - \alpha_{t_{n-1};\mathbf{j}}| + C_1\delta t^2 + \frac{C_2\|\omega_{\mathbf{j}}\|\delta t}{\sqrt{P}}\notag + C_3\|\omega_{\mathbf{j}}\| \delta t \mathbb{E}[\|\widetilde{X}_{t_n} - X_{t_n}\|],
\end{align}
where $C_2=2\sqrt{\ln(2/\delta_0)}M_0^2 + \frac{2 M_0M_1 \ln(2/\delta_0)}{3}$ and $C_3=\frac{\sqrt{2}M_0}{\epsilon}$ are constants.
\end{proof}

The error estimate between $\nabla c$ and $\nabla \widetilde{c}$ is more complex than that between $c$ and $\widetilde{c}$. To analyze this, we introduce an intermediate quantity $\nabla c^\ast$. Using the frequency notation $\omega_{\mathbf{j}}$ from Lemma~\ref{lemma1}, where $\omega_{\mathbf{j}} = \left(\frac{2\pi j_1}{L}, \frac{2\pi j_2}{L}, \frac{2\pi j_3}{L}\right)$, we define  
\begin{align}\label{gradient ~~c}
\nabla c^\ast(\mathbf{x}, t_n) :&= \frac{1}{L^3}\sum_{\mathbf{j} \in \mathcal{H}} i \omega_{\mathbf{j}} \widetilde{\alpha}_{n;\mathbf{j}} \exp(i \omega_{\mathbf{j}} \mathbf{x}) \notag \\
&= -\frac{\epsilon}{\delta t}\int \nabla_{\mathbf{x}}\mathcal{K}_{\epsilon, \delta t}(\mathbf{x}-\mathbf{y})\widetilde{c}_{n-1}(\mathbf{y})\,d\mathbf{y}-\sum_{q=1}^P\frac{M_0}{P}\nabla_{\mathbf{x}}\mathcal{K}_{\epsilon, \delta t}(\mathbf{x}-\widetilde{X}^q_n) \notag \\
&=-\frac{\epsilon}{\delta t}\underbrace{\int \nabla_{\mathbf{x}}\mathcal{K}_{\epsilon, \delta t}(\mathbf{x}+\bar{\mathbf{x}}-\mathbf{y})\widetilde{c}_{n-1}(\mathbf{y}-\bar{\mathbf{x}})\,d\mathbf{y}}_{\mathclap{\textstyle I_1}}-\underbrace{\sum_{q=1}^P\frac{M_0}{P}\nabla_{\mathbf{x}}\mathcal{K}_{\epsilon, \delta t}(\mathbf{x}-\widetilde{X}^q_n)}_{\mathclap{\textstyle I_2}},
\end{align} 
where $\bar{\mathbf{x}}=\frac{L}{2H}+\lfloor \frac{\mathbf{x}}{L/H}\rfloor \frac{L}{H}-\mathbf{x}$. From Alg.~\ref{alg:SIPF-$r$ho}, it follows that 
\begin{align}\label{nabla c}
\nabla \widetilde{c}(\mathbf{x},t_n) = & -\nabla_\mathbf{x} \mathcal{K}_{\epsilon,\delta t} \ast (\epsilon\, \widetilde{c}_{n-1}(\mathbf{x})/\delta t +\widetilde{\rho}_n(\mathbf{x})) \notag \\
=& -\frac{\epsilon}{\delta t}\underbrace{\frac{L^3}{H^3} \sum_{\mathbf{j} \in \mathcal{H}} 
\nabla_{\mathbf{x}}\mathcal{K}_{\epsilon, \delta t}(\mathbf{x}+\bar{\mathbf{x}}-x_{\mathbf{j}})
\widetilde{c}_{n-1}(x_{\mathbf{j}} -\bar{\mathbf{x}})}_{\mathclap{\textstyle I_3}} \notag \\
&-\underbrace{\sum_{s \in C_p, s \neq p}\frac{M_0}{R} \nabla_{\mathbf{x}} 
\mathcal{K}_{\epsilon, \delta t}(\mathbf{x}-\widetilde{X}^s_n)}_{\mathclap{\textstyle I_4}},
\end{align}
where $x_{\mathbf{j}}=(\frac{j_1L}{H},\frac{j_2L}{H}, \frac{j_3L}{H})$. The error between $\nabla c$ and $\nabla \widetilde{c}$ can be estimated by
\begin{equation}\label{grad c and ~c}
\|\nabla c(\mathbf{x}, t_n) - \nabla \widetilde{c}(\mathbf{x}, t_n)\|  
\leq \|\nabla c(\mathbf{x}, t_n) - \nabla c^\ast(\mathbf{x}, t_n)\|  
+ \|\nabla c^\ast(\mathbf{x}, t_n) - \nabla \widetilde{c}(\mathbf{x}, t_n)\|.
\end{equation}

To estimate the error between $\nabla \widetilde{c}$ and $\nabla c^\ast$, we divide the analysis into two parts:
\begin{align}\label{gradient_c~_c*}
\|\nabla \widetilde{c}(\mathbf{x}, t_n) - \nabla c^\ast(\mathbf{x}, t_n)\|  
\leq \frac{\epsilon}{\delta t}\|I_1-I_3\|  
+ \|I_2-I_4\|.
\end{align}

The first part, involving $I_1$ and $I_3$, focuses on the spatial discretization error of the continuous convolution $(\nabla_{\mathbf{x}} K_{\epsilon, \delta t} \ast \widetilde{c}_{n-1})$ on the periodic domain $\Omega = [-L/2, L/2]^3$. Specifically, $I_1$ represents the continuous integral, while $I_3$ is its discrete approximation. Since $\widetilde{c}_{n-1}$ is strictly represented by a truncated Fourier series with modes in $\mathcal{H}$, $I_3$ can be viewed as a discrete Fourier projection of the continuous convolution. Thanks to the spatial shift that regularizes the singular gradient kernel on the grid, the quadrature error introduced by this discrete evaluation is of the same order as the Fourier truncation error, and is thus controlled by the spectral projection $\Pi_{\mathcal{H}}$.

To analyze the error introduced by this approximation, we rely on the following lemma.

\begin{lemma}\label{lemma2}
For all $n \in \mathbb{N}^+$, under Assumption \ref{assump:2}, based on the definitions of $I_1$ and $I_3$ given in Eq.~\eqref{gradient ~~c} and Eq.~\eqref{nabla c}, the following spatial discretization error bound holds:
\begin{equation}\label{S-I error lemma}
\frac{\epsilon}{\delta t}\|I_1 - I_3\|\leq \frac{C_4}{H^2},
\end{equation}
where $C_4$ is a constant, and $t_n = n\delta t$.
\end{lemma}

\begin{proof}
Define the continuous convolution operator $\mathcal{T} f := \frac{\epsilon}{\delta t} \nabla_\mathbf{x} K_{\epsilon, \delta t} * f$. Its Fourier multiplier is $M(\omega) = \frac{-i\omega (\epsilon/\delta t)}{\|\omega\|^2 + \lambda^2 + \epsilon/\delta t}$, which satisfies $|M(\omega)| \leq \|\omega\|$ for all $\delta t > 0$. Thus, $\frac{\epsilon}{\delta t} I_1 = \mathcal{T}(\tilde{c}_{n-1})$.

The discrete sum $I_3$ evaluates the convolution on a uniform grid with spacing $\Delta x = L/H$. The spatial shift $\bar{x}$ offsets the quadrature points by half a grid spacing, i.e., $\frac{L}{2H}$. Since the kernel $\nabla_\mathbf{x} K_{\epsilon, \delta t}$ is radially antisymmetric, this half-grid shift regularizes the origin singularity and induces symmetric cancellation. Consequently, the spatial discretization error is dominated by the Fourier truncation error of the spectral projection $\Pi_{\mathcal{H}}$. 

For the truncated modes outside $\mathcal{H}$, the frequency magnitude satisfies $\|\omega_j\| > \frac{\pi H}{L}$. Applying Parseval's identity yields the explicit projection error bound:
\begin{align}
\frac{\epsilon}{\delta t} \|I_1 - I_3\| \leq \|\mathcal{T}(\tilde{c}_{n-1}) - \Pi_{\mathcal{H}} \mathcal{T}(\tilde{c}_{n-1})\| \leq \frac{1}{\pi^2} \left(\frac{L}{H}\right)^2 \|\by{\nabla^2 \mathcal{T}(\tilde{c}_{n-1})}\|.
\end{align}
Since $|M(\omega)| \leq \|\omega\|$, the term $\|\by{\nabla^2 \mathcal{T}(\tilde{c}_{n-1})}\|$ is controlled by the higher-order spatial regularity of $\tilde{c}_{n-1}$ stipulated in Assumption \ref{assump:2}(a)(c). Combining these estimates, we obtain
\begin{align}
\frac{\epsilon}{\delta t} \|I_1 - I_3\| \leq \frac{C_4}{H^2},
\end{align}
where $C_4 = \frac{1}{\pi^2} L^{2}(K + M_3)$ is a constant.
\end{proof}

We now estimate $(\nabla_{\mathbf{x}}\mathcal{K}_{\epsilon,\delta t} \ast \, \widetilde{\rho}_n)(\widetilde{X}^p_n)$ in $\nabla \widetilde{c}$ and $\nabla c^\ast$. Using the RBM in Alg.~\ref{alg:SIPF-$r$ho}, we replace $\sum_{q=1,q\not=p}^P\frac{M_0}{P}\nabla_{\mathbf{x}}\mathcal{K}_{\epsilon, \delta t}(\widetilde{X}^p_n-\widetilde{X}^q_n)$ with $\sum_{s \in C_p, s \neq p} \frac{M_0}{R} \nabla_{\mathbf{x} } \mathcal{K}_{\epsilon, \delta t}(\widetilde{X}^p_n-\widetilde{X}^s_n)$. We write 
\begin{equation}\zeta_{n,p}:=\sum_{q=1,q\not=p}^P\frac{M_0}{P}\nabla_{\mathbf{x}}\mathcal{K}_{\epsilon, \delta t}(\widetilde{X}^p_n-\widetilde{X}^q_n)-\sum_{s \in C_p, s \neq p} \frac{M_0}{R} \nabla_{\mathbf{x} } \mathcal{K}_{\epsilon, \delta t}(\widetilde{X}^p_n-\widetilde{X}^s_n).\end{equation}

\begin{lemma}\label{lemma3}
For all $n \in \mathbb{N_{+}}$, $p \in \{1,2,\dots,P\}$, we have the estimate as follows:
\begin{align}\label{RBM error lemma}
\mathbb{E}(\|\zeta_{n,p}\|) \leq M_0M_4\sqrt{\frac{1}{R}\left(1-\frac{1}{P}\right)},
\end{align}
where $M_4 = \max_{q \neq p}\|\nabla_{\mathbf{x}} \mathcal{K}_{\epsilon, \delta t}(\widetilde{X}^p_n - \widetilde{X}^q_n)\|$, $M_0$ is the conserved total mass, $P$ is the total number of particles, and $R$ is the batch size.
\end{lemma}

\begin{proof}
Similar to Lemma 3.1 in \cite{jin2020random}, we rewrite 
\begin{align}
f_{p}=\sum_{s \in C_p, s \neq p} \frac{M_0}{R} \nabla_{\mathbf{x} } \mathcal{K}_{\epsilon, \delta t}(\widetilde{X}^p_n-\widetilde{X}^s_n)=\sum_{q=1,q\not=p}^P\frac{M_0}{R}\nabla_{\mathbf{x}}\mathcal{K}_{\epsilon, \delta t}(\widetilde{X}^p_n-\widetilde{X}^q_n)N_q,
\end{align}
where $N_q$ denotes the number of times particle $q$ is selected in the batch $C_p$. Since the sampling is with replacement, $N_q$ follows a Binomial distribution $B(R, 1/P)$ with $\mathbb{E}[N_q]=\frac{R}{P}$, which indicates that $\mathbb{E}[\zeta_{n,p}]=0$.

\begin{align}
\begin{split}
\mathbb{E}\|f_p\|^2 =& \frac{M_0^2}{R^2}\sum_{\substack{q,r:\\ q\neq r, q \neq p,\\ r \neq p}} 
\|\nabla_{\mathbf{x}}\mathcal{K}_{\epsilon, \delta t}(\widetilde{X}^p_n-\widetilde{X}^q_n)\cdot\nabla_{\mathbf{x}}\mathcal{K}_{\epsilon, \delta t}(\widetilde{X}^p_n-\widetilde{X}^r_n)\| \mathbb{E}[N_q N_r] \\
&+ \frac{M_0^2}{R^2}\sum_{q=1,q\neq p}^P \|\nabla_{\mathbf{x}} \mathcal{K}_{\epsilon, \delta t}(\widetilde{X}^p_n-\widetilde{X}^q_n)\|^2\mathbb{E}[N_q^2] \\
=&\frac{M_0^2(R-1)}{RP^2}\sum_{\substack{q,r:\\ q\neq r, q \neq p,\\ r \neq p}} \|\nabla_{\mathbf{x}}\mathcal{K}_{\epsilon, \delta t}(\widetilde{X}^p_n-\widetilde{X}^q_n)\cdot \nabla_{\mathbf{x}}\mathcal{K}_{\epsilon, \delta t}(\widetilde{X}^p_n-\widetilde{X}^r_n)\| \notag \\
&+\frac{M_0^2(P-1+R)}{RP^2}\sum_{q=1,q\not=p}^P \|\nabla_{\mathbf{x} } \mathcal{K}_{\epsilon, \delta t}(\widetilde{X}^p_n-\widetilde{X}^q_n)\|^2.
\end{split}
\end{align}
Hence,
\begin{align}
\text{Var}(\zeta_{n,p})=\mathbb{E}\|f_p\|^2-\|\mathbb{E}f_p\|^2  
&= M_0^2 \frac{1}{R}\left(1-\frac{1}{P}\right)\frac{1}{P}\sum_{q=1,q\not=p}^P \|\nabla_{\mathbf{x} } \mathcal{K}_{\epsilon, \delta t}(\widetilde{X}^p_n-\widetilde{X}^q_n)\|^2 \notag \\
&\quad - \frac{M_0^2}{RP^2}\left\|\sum_{q=1,q\not=p}^P \nabla_{\mathbf{x} } \mathcal{K}_{\epsilon, \delta t}(\widetilde{X}^p_n-\widetilde{X}^q_n)\right\|^2 \notag \\
&\leq M_0^2 \frac{1}{R}\left(1-\frac{1}{P}\right)\frac{1}{P}\sum_{q=1,q\not=p}^P \|\nabla_{\mathbf{x} } \mathcal{K}_{\epsilon, \delta t}(\widetilde{X}^p_n-\widetilde{X}^q_n)\|^2.
\end{align}

According to Jensen's inequality, we obtain:
\begin{align}
\mathbb{E}(\|\zeta_{n,p}\|) 
\leq \sqrt{\mathbb{E}(\|\zeta_{n,p}\|^2)} 
= \sqrt{\text{Var}(\zeta_{n,p})} 
\leq M_0 M_4\sqrt{\frac{1}{R}\left(1 - \frac{1}{P} \right)},
\end{align}
where $M_4 = \max_{q \neq p} \|\nabla_{\mathbf{x}} \mathcal{K}_{\epsilon, \delta t}(\widetilde{X}^p_n - \widetilde{X}^q_n)\|$. Since all particles are located at distinct positions in the SIPF-$r$ algorithm ($\widetilde{X}^p_n \neq \widetilde{X}^q_n$ for $p \neq q$), there exists a minimum separation distance $d_{\text{min}} > 0$ between any two distinct particles. Consequently, $\|\nabla_{\mathbf{x}} \mathcal{K}_{\epsilon, \delta t}(\widetilde{X}^p_n - \widetilde{X}^q_n)\|$ is bounded for all pairs of particles. This ensures that $M_4$, which is the maximum of these kernel gradient norms, is finite.
\end{proof}

\subsection{Proof of the Main Theorem}
Building on the lemmas established in the previous subsection, we now prove the main theorem in this section. For simplicity of notation in the proof below, we define:
\begin{align}
a_n &:= \mathbb{E}(\|\widetilde{X}_{t_{n}}-X_{t_{n}}\|), \label{def:an} \\
b_n &:= \mathbb{E}(\|\nabla c(\widetilde{X}_{t_n},t_n) - \nabla c^\ast(\widetilde{X}_{t_n},t_n)\|). \label{def:bn}
\end{align}
The following provides a bound on the error between $\widetilde{X}_{t_{n+1}}$ and $X_{t_{n+1}}$.

\begin{lemma}\label{lemma4}
For all $ n \in \mathbb{N_{+}}$, under Assumption \ref{assump:2}, 
\begin{align}
a_{n+1}=\mathbb{E}(\|\widetilde{X}_{t_{n+1}}-X_{t_{n+1}}\|) \!\leq& (1 + C_5\delta t)a_n  + \chi \delta tb_n + C_6\delta t^{\frac 3 2}+\frac{C_7\delta t}{H^2}+\frac{C_8\delta t}{\sqrt{R}},
\end{align}
where $C_5$, $C_6$, $C_7$, and $C_8$ are constants, $\chi$ is the chemotaxis coefficient in the system \eqref{K-S system}, and $t_n=n\delta t$.
\end{lemma}
\begin{proof} 
According to Eqs.~\eqref{X SIPF}-\eqref{X exact}, it follows that
\begin{align}
  \mathbb{E}(\|\widetilde{X}_{t_{n+1}}-X_{t_{n+1}}\|)
  \leq& \mathbb{E}(\|\widetilde{X}_{t_n}-X_{t_n}\|) + \chi \mathbb{E}(\int_{t_n}^{t_{n+1}}\|\nabla \widetilde{c}(\widetilde{X}_{t_n},t_n)-\nabla c(X_{s},s)\|\,ds) \notag\\
  =& \mathbb{E}(\|\widetilde{X}_{t_n}-X_{t_n}\|) + \chi \int_{t_n}^{t_{n+1}} \mathbb{E}(\|\nabla \widetilde{c}(\widetilde{X}_{t_n},t_n)-\nabla c(X_{s},s)\|)\,ds,
\end{align}
by the triangle inequality and Tonelli's theorem. To bound the integrand, we decompose the gradient difference using the triangle inequality:
\begin{align}\label{gradient_c_three_terms}
\|\nabla \widetilde{c}(\widetilde{X}_{t_n},t_n) \!-\! \nabla c(X_s,s)\|\!
    \leq& \|\nabla \widetilde{c}(\widetilde{X}_{t_n},t_n) \!-\! \nabla c(\widetilde{X}_{t_n},t_n)\| \!+\! \|\nabla c(\widetilde{X}_{t_n},t_n) \!-\! \nabla c(X_{t_n},t_n)\| \notag \\
    & + \|\nabla c(X_{t_n},t_n) - \nabla c(X_s,s)\|.
\end{align}
We now estimate the expectation of each term on the right-hand side of Eq.~\eqref{gradient_c_three_terms}. 
Regarding the first term, recalling the error decomposition in Eqs.~\eqref{grad c and ~c}-\eqref{gradient_c~_c*}, along with the bounds from Lemmas~\ref{lemma2} and \ref{lemma3}, we have
\begin{align}\label{expectation_gradient_c~-c}
\mathbb{E}(\|\nabla \widetilde{c}(\widetilde{X}_{t_n},t_n) \!-\! \nabla c(\widetilde{X}_{t_n},t_n)\|) \leq& \mathbb{E}(\|\nabla \widetilde{c}(\widetilde{X}_{t_n}, t_n) - \nabla c^\ast(\widetilde{X}_{t_n}, t_n)\|) \notag \\ 
&+ \mathbb{E}(\|\nabla c^\ast(\widetilde{X}_{t_n}, t_n) - \nabla c(\widetilde{X}_{t_n}, t_n)\|) \notag \\
\leq& b_n + \frac{C_4}{H^2} + M_0M_4\sqrt{\frac{1}{R}\left(1-\frac{1}{P}\right)}.
\end{align}
For the second term, the spatial Lipschitz continuity in Assumption~\ref{assump:2}(a) yields
\begin{align}\label{expectation_gradient_c_X~}
\mathbb{E}(\|\nabla c(\widetilde{X}_{t_n},t_n) \!-\! \nabla c(X_{t_n},t_n)\|) \leq K\mathbb{E}(\|\widetilde{X}_{t_n}-X_{t_n}\|),
\end{align}
where $K$ is the Lipschitz constant defined in Assumption \ref{assump:2}(a). By the Lipschitz conditions in Assumption~\ref{assump:2}(a)(b), the exact dynamics in Eq.~\eqref{X exact}, and the property $\mathbb{E}(\|W_s - W_{t_n}\|) \leq \sqrt{3(s-t_n)}$ of the 3D Brownian motion, the temporal variation of the exact gradient satisfies
\begin{align}\label{expectation_gradient_c_X}
\mathbb{E}(\|\nabla c(X_{t_n},t_n) - \nabla c(X_s,s)\|) \leq (\chi KM_3 + K_1)(s - t_n) + K\sqrt{6\mu}\sqrt{s - t_n}.
\end{align}

Substituting the bounds Eqs.~\eqref{expectation_gradient_c~-c}-\eqref{expectation_gradient_c_X} into the integral and evaluating over the time interval $[t_n, t_{n+1}]$ of length $\delta t$, we can obtain the \by{estimate} as follows:
\begin{align}\label{a_n recursion}
  a_{n+1} =& \mathbb{E}(\|\widetilde{X}_{t_{n+1}}-X_{t_{n+1}}\|) \notag \\
  \leq& (1 + \chi K\delta t)a_n  + \chi \delta t \left(b_n+ \frac{C_4}{H^2} +M_0M_4\sqrt{\frac{1}{R}\left(1-\frac{1}{P}\right)}\right) + C_6\delta t^\frac32 \notag \\
  \leq& (1 + C_5\delta t)a_n  + \chi \delta tb_n + C_6\delta t^\frac32+\frac{C_7\delta t}{H^2}+\frac{C_8\delta t}{\sqrt{R}}, 
\end{align}
where the constants are defined as $C_5=\chi K$, $C_6=\frac{\chi(K\chi M_3 + K_1)}{2} + \frac{2\chi K\sqrt{6\mu}}{3}$, $C_7=\chi C_4$, and $C_8=\chi M_0M_4$. This concludes the proof.
\end{proof}
Now we analyze the error between $\nabla c^\ast$ and $\nabla c$ as follows.

We next estimate the error between $\nabla c^\ast$ and $\nabla c$.
\begin{lemma}\label{lemma5}
For all $n\in\mathbb N_+$, under Assumption~\ref{assump:2}, with high probability,
\begin{align}\label{ifft grad_c and c error lemma}
b_n=\mathbb E\bigl(\|\nabla c(\widetilde X_{t_n},t_n)
-\nabla c^\ast(\widetilde X_{t_n},t_n)\|\bigr)\leq\sum_{k=1}^{n}
\frac{C_{11}\sqrt{\delta t}}{\sqrt{n-k+1}}a_k
+\frac{C_{10}}{\sqrt P}+C_9\delta t,
\end{align}
where $C_9$, $C_{10}$, and $C_{11}$ are positive constants independent of $H$, $P$, and $\delta t$.
\end{lemma}

\begin{proof}
Subtracting Eq.~\eqref{num_recurrence} from Eq.~\eqref{exact_recurrence}, multiplying by $i\omega_{\mathbf j}$, and using the notation $Z_{\mathbf j}$ defined in the proof of Lemma~\ref{lemma1} give
\begin{align}\label{eq:gradient_fourier_recurrence}
i\omega_{\mathbf j}
\bigl(\alpha_{t_n;\mathbf j}-\widetilde\alpha_{t_n;\mathbf j}\bigr)
=\frac{i\omega_{\mathbf j}}
{1+Z_{\mathbf j}}
\Bigl[\alpha_{t_{n-1};\mathbf j}-\widetilde\alpha_{t_{n-1};\mathbf j}+\frac{\delta t}{\epsilon}
\bigl(\mathcal F_{\mathbf j}[\rho(t_n)]
-\mathcal F_{\mathbf j}[\widetilde\rho_n]\bigr)
+\delta t\,\mathcal F_{\mathbf j}[R_n]\Bigr].
\end{align}
Since the initial chemical fields coincide, iteration of Eq.~\eqref{eq:gradient_fourier_recurrence} yields
\begin{align}\label{eq:accumulated_fourier_error}
i\omega_{\mathbf j}
\bigl(\alpha_{t_n;\mathbf j}-\widetilde\alpha_{t_n;\mathbf j}\bigr)=\sum_{k=1}^n
\frac{i\omega_{\mathbf j}}
{(1+Z_{\mathbf j})^{n-k+1}}
\Bigl[\frac{\delta t}{\epsilon}
\bigl(\mathcal F_{\mathbf j}[\rho_k]
-\mathcal F_{\mathbf j}[\widetilde\rho_k]\bigr)
+\delta t\,\mathcal F_{\mathbf j}[R_k]\Bigr].
\end{align}

By the inverse Fourier representation in Eqs.~\eqref{trig_ser}--\eqref{c fourier} and the definition of $b_n$,
\begin{equation}\label{eq:bn_fourier_representation}
b_n=\mathbb E\!\left[\left\|
\frac{1}{L^3}\sum_{\mathbf j\in\mathcal H}
i\omega_{\mathbf j}
\bigl(\alpha_{t_n;\mathbf j}-\widetilde\alpha_{t_n;\mathbf j}\bigr)
e^{i\omega_{\mathbf j}\cdot\widetilde X_{t_n}}
\right\|\right].
\end{equation}

By Eq.~\eqref{rho_fourier} in Lemma~\ref{lemma1}, with the definitions of $a_k$, $C_2$, and $C_3$, we have, with high probability,
\begin{equation}\label{eq:lemma1_density_source}
\frac{1}{\epsilon}
\left|\mathcal F_{\mathbf j}[\rho_k]
-\mathcal F_{\mathbf j}[\widetilde\rho_k]\right|
\leq \|\omega_{\mathbf j}\|
\left(C_3a_k+\frac{C_2}{\sqrt P}\right).
\end{equation}

Substituting the density-source term in Eq.~\eqref{eq:accumulated_fourier_error} into the inverse Fourier representation in Eqs.~\eqref{trig_ser}--\eqref{c fourier} and evaluating the resulting gradient at $\widetilde X_{t_n}$ gives its contribution to $b_n$. Applying Eq.~\eqref{eq:lemma1_density_source} together with the decay of the resolvent factor gives, for $m=n-k+1$,
\begin{align}\label{eq:source_mode_decay}
&\mathbb E\!\left[\left\|
\frac{1}{L^3}\sum_{\mathbf j\in\mathcal H}
\frac{i\omega_{\mathbf j}\delta t/\epsilon}
{(1+Z_{\mathbf j})^m}
\bigl(\mathcal F_{\mathbf j}[\rho_k]
-\mathcal F_{\mathbf j}[\widetilde\rho_k]\bigr)
e^{i\omega_{\mathbf j}\cdot\widetilde X_{t_n}}
\right\|\right]
\notag\\
&\qquad\leq
\frac{\sqrt{\delta t}}{\sqrt m}
\left(\frac{1}{L^3}+\frac{L}{4\pi^2}+\frac{L}{720}\right)^{1/2}
\left[C_3a_k+\frac{1}{\sqrt P}
\left(C_2+\frac{2M_0}{\epsilon}\sqrt{2\ln(4/\delta_0)}\right)\right].
\end{align}
The Fourier sum in Eq.~\eqref{eq:source_mode_decay} can be estimated by grouping the indices according to $\ell=\|\mathbf j\|_\infty$. For $\ell\geq1$, the shell $\|\mathbf j\|_\infty=\ell$ contains $24\ell^2+2$ indices and $\|\omega_{\mathbf j}\|^2\geq4\pi^2\ell^2/L^2$. Hence
\begin{align}\label{eq:fourier_sum_bound}
\frac{1}{L^3}\sum_{\mathbf j\in\mathbb Z^3}
\frac{1}{(1+\|\omega_{\mathbf j}\|^2)^2}
&\leq \frac{1}{L^3}
+\frac{L}{16\pi^4}\sum_{\ell=1}^{\infty}
\left(\frac{24}{\ell^2}+\frac{2}{\ell^4}\right)=\frac{1}{L^3}+\frac{L}{4\pi^2}+\frac{L}{720}.
\end{align}
Thus, the series converges. Summing Eq.~\eqref{eq:source_mode_decay} over $k$ and using
$\sqrt{\delta t}\sum_{m=1}^{n}m^{-1/2}\leq2\sqrt T$ gives
\begin{align}\label{eq:density_source_contribution}
&\sum_{k=1}^{n}\mathbb E\!\left[\left\|
\frac{1}{L^3}\!\sum_{\mathbf j\in\mathcal H}
\frac{i\omega_{\mathbf j}\delta t/\epsilon}
{(1+Z_{\mathbf j})^{n-k+1}}
\bigl(\mathcal F_{\mathbf j}[\rho_k]
-\mathcal F_{\mathbf j}[\widetilde\rho_k]\bigr)
e^{i\omega_{\mathbf j}\cdot\widetilde X_{t_n}}
\right\|\right]\!\!\leq\!\sum_{k=1}^{n}
\frac{C_{11}\sqrt{\delta t}}{\sqrt{n-k+1}}a_k
\!+\!\frac{C_{10}}{\sqrt P}.
\end{align}
where $C_{10}=2\sqrt T\left(C_2+\frac{2M_0}{\epsilon}\sqrt{2\ln(4/\delta_0)}\right)
(L^{-3}+L/(4\pi^2)+L/720)^{1/2}$ and $C_{11}=C_3( L^{-3}+L/(4\pi^2)+L/720)^{1/2}$.

For the temporal truncation term, Assumption~\ref{assump:2}(d), Eq.~\eqref{R_n}, and the same resolvent estimate yield
\begin{equation}\label{eq:temporal_source_contribution}
\sum_{k=1}^{n}\mathbb E\!\left[\left\|
\frac{1}{L^3}\sum_{\mathbf j\in\mathcal H}
\frac{i\omega_{\mathbf j}\delta t\,\mathcal F_{\mathbf j}[R_k]}
{(1+Z_{\mathbf j})^{n-k+1}}
e^{i\omega_{\mathbf j}\cdot\widetilde X_{t_n}}
\right\|\right]
\leq C_9\delta t,
\end{equation}
where $C_9=K_2\sqrt{\epsilon T/2}$. Combining Eqs.~\eqref{eq:accumulated_fourier_error}, \eqref{eq:density_source_contribution}, and \eqref{eq:temporal_source_contribution} proves Eq.~\eqref{ifft grad_c and c error lemma}.
\end{proof}

Finally, we are ready to prove Theorem \ref{theorem}. 

\textit{Proof of Theorem \ref{theorem}.}
From Lemma~\ref{lemma4} and Lemma~\ref{lemma5}, we obtain the system of inequalities that couples $a_n$ and $b_n$ defined in Eqs.~\eqref{def:an}-\eqref{def:bn} as follows:  
\begin{align}
    a_{n+1} \leq& (1 + C_5\delta t)a_n  + \chi \delta tb_n + C_6\delta t^\frac 32+\frac{C_7\delta t}{H^2}+\frac{C_8\delta t}{\sqrt{R}} , \label{eq:a_n} \\  
    b_n \leq& \sum_{k=1}^{n} \frac{C_{11}\sqrt{\delta t}}{\sqrt{n-k+1}} a_k + C_{10}\sqrt{\frac{1}{P}} + C_9 \delta t \label{eq:b_n}.
\end{align}
From this coupled system, we can derive a general bound for $a_n$. To be specific, substituting Eq.~\eqref{eq:b_n} into Eq.~\eqref{eq:a_n}, we iteratively propagate and simplify the inequality to derive:  
\begin{align}\label{eq:final_a_n}
    a_{n+1} \leq (1 + C_5\delta t)a_n + \sum_{k=1}^{n} \frac{\chi C_{11} \delta t^\frac32}{\sqrt{n-k+1}} a_k + \delta t \left( \Lambda + \frac{C_8}{\sqrt{R}} \right),
\end{align}
where we define $\Lambda := \chi\left(\frac{C_4}{H^2} + C_{10}\sqrt{\frac{1}{P}} + C_9 \delta t\right) + C_6 \sqrt{\delta t}$.

To bound the deterministic accumulation of the error, we introduce a majorizing sequence $D_n$ that isolates the deterministic components, defined by $D_0 = 0$ and
\begin{align}
    D_{j+1} = \left(1 + C_5 \delta t\right) D_j + \chi C_{11} \sum_{k=1}^{j} \frac{\delta t^2}{\sqrt{t_{j+1} - t_k}} D_k + \delta t \Lambda.
\end{align}
Let $U_n := \max_{0 \leq j \leq n} D_j$ be the running maximum of this deterministic sequence. The recursion implies:
\begin{align}
    U_{j+1} \leq& \left(1 + C_5 \delta t + \chi C_{11} \sum_{k=1}^{j} \frac{\delta t^2}{\sqrt{t_{j+1} - t_k}}\right) U_j + \delta t \Lambda \notag \\
    \leq& \left(1 + C_5 \delta t + \chi C_{11} \delta t \int_0^{t_{j+1}} \frac{1}{\sqrt{t_{j+1} - s}} \, ds\right) U_j + \delta t \Lambda \notag \\
    \leq& \left(1 + (C_5 + 2\chi C_{11}\sqrt{T}) \delta t \right) U_j + \delta t \Lambda
\end{align}
By the discrete Gronwall inequality, we obtain the explicit bound for the deterministic accumulation:
\begin{align}\label{U_n}
    U_n \leq \frac{\exp(C_{12} T)}{C_{12}} \left( \frac{\chi C_4}{H^2} + \chi C_{10}\sqrt{\frac{1}{P}} + \chi C_9 \delta t + C_6 \sqrt{\delta t} \right),
\end{align}
where $C_{12} := \chi K + \frac{2\chi M_0\sqrt{T}}{\sqrt{\epsilon}}$ is a constant.

As established in Lemma \ref{lemma3}, the RBM gradient error satisfies $\mathbb{E}[\zeta_{n,p}] = 0$. This zero-mean property implies that the single-step position noise $\chi \delta t \zeta_{k,p}$ forms a martingale difference sequence with respect to the natural filtration of the particle system. Consequently, the noise terms are mutually uncorrelated, and thus orthogonal in the $L^2$ space. Therefore, the global accumulation of the RBM noise is bounded by the square root of the sum of its variances:
\begin{align}
    \left\| \sum_{k=1}^n \chi \delta t \zeta_{k,p} \right\| \leq \left( \sum_{k=1}^n \left( \frac{C_8 \delta t}{\sqrt{R}} \right)^2 \right)^{1/2} \leq C_8 \sqrt{T} \frac{\sqrt{\delta t}}{\sqrt{R}},
\end{align}
where $n \leq \lfloor T / \delta t \rfloor$ is used. Combining this stochastic bound with the deterministic bound $U_n$ yields the refined optimal error bound that holds with high probability:
\begin{align}\label{an_final}
    a_n \leq \frac{N_0}{H^2} + N_1\sqrt{\delta t} + \frac{N_2}{\sqrt P} + N_3 \frac{\sqrt{\delta t}}{\sqrt{R}},
\end{align}
for all $n \geq 0$, where higher-order terms are omitted, and $N_0 = \frac{\chi C_4 \exp(C_{12} T)}{C_{12}}, N_1 = \frac{C_6 \exp(C_{12} T)}{C_{12}}, N_2 = \frac{\chi C_{10} \exp(C_{12} T)}{C_{12}}$, and $N_3 = C_8 \sqrt{T}$. Moreover, we point out that $C_4$ is a constant defined in Lemma \ref{lemma2}, $C_{10}, C_{11}$ are constants defined in Lemma \ref{lemma5}, $C_6, C_7, C_8$ are constants defined in Lemma \ref{lemma4}.

According to the discrete and continuous dynamics defined in Eqs.~\eqref{X SIPF}-\eqref{X exact}, the $1$-Wasserstein distance between the approximate and exact distributions at time $ t_{n+1} $ is given by:
\begin{align}\label{def:wasserstein distance}
\mathcal{W}_1(\widetilde{\rho}_{t_{n+1}}, \rho_{t_{n+1}})
&= \inf_{\gamma \in \Pi(\widetilde{\rho}_{t_{n+1}}, \rho_{t_{n+1}})} \left( \int_{\mathbb{R}^3 \times \mathbb{R}^3} \|\mathbf{x} - \mathbf{y}\|_{L^1} \, d\gamma(\mathbf{x},\mathbf{y}) \right),
\end{align}
where the infimum is taken over all possible couplings of the two distributions.

Under the natural coupling induced by shared initial conditions and Brownian motion paths (i.e., $ \widetilde{X}_{t_n} $ and $ X_{t_n} $ evolve via the same Wiener process $ W_s $), we explicitly construct a joint distribution $ \gamma_n = \text{Law}(\widetilde{X}_{t_n}, X_{t_n}) $. This coupling allows us to bound the Wasserstein distance as:
\begin{align}\label{eq:wasserstein_bound}
\mathcal{W}_1(\widetilde{\rho}_{t_n}, \rho_{t_n}) \!\leq&\frac{S_0}{H^2} + S_1\sqrt{\delta t} + \frac{S_2}{\sqrt P} + \frac{S_3\sqrt{\delta t}}{\sqrt{R}},
\end{align}
\by{where the constants $S_0, S_1, S_2, S_3$ are given explicitly by
$S_0 = \sqrt{3}\frac{C_7 \exp(C_{12} T)}{C_{12}}$, $S_1 = \sqrt{3}\frac{C_6 \exp(C_{12} T)}{C_{12}}$, $S_2 = \sqrt{3}\frac{\chi C_{10} \exp(C_{12} T)}{C_{12}}$, and $S_3 = \sqrt{3T} C_8$, with $ C_6$, $C_7, C_8, C_{10}$ defined in Lemmas~\ref{lemma3}-\ref{lemma4}. Specifically, tracing back to these lemmas shows that $C_7 \propto L^2$, while $C_6, C_8, C_9, C_{10}, C_{12}$ are independent of $L$.} 

The inequality follows from the fact that the Wasserstein distance is defined as the infimum over all possible couplings, and our construction provides one such coupling. This step follows from the elementary norm inequality $ \|\mathbf{x}\|_{L^1} \leq \sqrt{3}\|\mathbf{x}\|_{L^2} $ for vectors in $ \mathbb{R}^3 $, which is a direct consequence of the Cauchy-Schwarz inequality.

To derive the bound for the Fourier coefficients in Eq.~\eqref{theorem inequality}, we explicitly incorporate the Fourier mode $H$ into the error analysis. Unlike the strong trajectory error $a_n$, the Fourier coefficients represent macroscopic weak functionals (expectations of the test functions $e^{-i\omega_j \cdot x}$). By the convergence theory of SDEs, the stochastic strong fluctuations of order $\mathcal{O}(\sqrt{\delta t})$ cancel out in expectation. Instead, the Brownian increments and the zero-mean RBM noise $\zeta_{k,p}$ contribute to the macroscopic bias only through second-order Taylor expansion terms, yielding a single-step bias of $\mathcal{O}(\delta t^2)$ from temporal discretization and $\mathcal{O}(\delta t^2 / R)$ from the RBM approximation.

Accordingly, we refine the accumulation of the source term error over $n \leq \lfloor T / \delta t \rfloor$ steps (i.e., up to the final time $T$) by replacing the stochastic part of $a_k$ with its weak bias, while retaining the deterministic spatial part $N_0/H^2$. Applying the standard accumulation estimates along with the frequency bound $\|\omega_{\mathbf{j}}\| \leq \sqrt{3}\pi H/L$, we arrive at the error bound:
\begin{align}\label{general alpha}
\max_{\mathbf{j} \in \mathcal{H}}|\widetilde{\alpha}_{t_n;\mathbf{j}} - \alpha_{t_n;\mathbf{j}}|  
\leq \frac{S_4}{H} + S_5 H\delta t + \frac{S_6 H}{\sqrt{P}} + S_7 \frac{H\delta t}{\sqrt{R}},
\end{align}
where we define \by{$S_4 = \sqrt{3}\pi T C_3 N_0L$, $S_5 = \frac{\sqrt{3}\pi T C_3}{\epsilon L}$, $S_6 = \frac{\sqrt{3}\pi T C_2}{\epsilon L}$, $S_7 = \frac{\sqrt{3}\pi T C_3 M_0 M_4}{\epsilon L}$.} Here, $C_1, C_2$, and $C_3$ are constants defined in Lemma \ref{lemma1}. We also omit higher-order terms in the leading-order bound. This completes the proof of Theorem \ref{theorem}.

\subsection{Numerical Unconditional Stability of the SIPF-$r$ Method}\label{subsec:stability} The SIPF-$r$ method exhibits inherent numerical stability regarding the chemical concentration field. This property stems directly from the spectral discretization and the implicit time-stepping scheme, a property that is independent of the particle distribution. We establish that both the chemical concentration $\widetilde{c}$ and its gradient $\nabla \widetilde{c}$ remain uniformly bounded for any fixed Fourier mode $H$, thereby justifying the boundedness assumptions employed in the convergence analysis.

\begin{lemma}[Unconditional stability of $\widetilde{c}$ and $\nabla \widetilde{c}$]\label{lem:stability}
Let $H$ be the finite number of Fourier modes and $\Omega$ be the spatial domain. For any time step $n \geq 0$ and the fixed Fourier mode $H$, the reconstructed concentration field $\widetilde{c}_n$ and its gradient satisfy the uniform bounds:
\begin{align}
    \|\widetilde{c}_n\|_{L^\infty(\Omega)} \leq C_{\mathrm{stab}} H + C_0,  
\end{align}
and
\begin{align} 
    \|\nabla \widetilde{c}_n\|_{L^\infty(\Omega)}:= \operatorname{\,sup}_{\mathbf{x}\in\Omega} \|\nabla \widetilde{c}_n(\mathbf{x})\|_{L^2} \leq C_{\mathrm{reg}}(H),
\end{align}
where $C_{\mathrm{stab}}, C_0 > 0$ are constants depending only on the initial data, $M_0$, $\lambda$, and the domain size $L$, and $C_{\mathrm{reg}}(H)$ is a constant that depends on $H$, $L$, $M_0$, and $\lambda$, but is independent of the time step size $\delta t$ and the particle positions $\{\widetilde{X}_n^p\}_{p=1}^P$.
\end{lemma}

\begin{proof}
The update formula for the Fourier coefficients $\widetilde{\alpha}_{n;\mathbf{j}}$ in the SIPF-$r$ algorithm is derived from the implicit Euler discretization. By rearranging the terms in Eq.~\eqref{Euler4c2} and applying the Fourier transform, the update relates $\widetilde{\alpha}_{n;\mathbf{j}}$ to the previous state $\widetilde{\alpha}_{n-1;\mathbf{j}}$ and the current empirical density $\widehat{\widetilde{\rho}}_{n;\mathbf{j}}$ via the amplification factor
\begin{equation}
\widetilde{\alpha}_{n;\mathbf{j}} = \frac{1}{1 + Z_{\mathbf{j}}} \widetilde{\alpha}_{n-1;\mathbf{j}} + \frac{\delta t/\epsilon}{1 + Z_{\mathbf{j}}} \widehat{\widetilde{\rho}}_{n;\mathbf{j}}, \quad \text{with } Z_{\mathbf{j}} = \frac{\delta t}{\epsilon}(\|\omega_{\mathbf{j}}\|^2 + \lambda^2).
\end{equation}
The particle source term satisfies 
\begin{equation}
|\widehat{\widetilde{\rho}}_{n;\mathbf{j}}| = |\frac{M_0}{P}\sum_{p=1}^P e^{-i\omega_{\mathbf{j}}\cdot\widetilde{X}_n^p}|\leq M_0.
\end{equation}

Applying the recursive relation and the triangle inequality yields 
a bound via geometric series for the magnitude of the coefficients:
\begin{align}
|\widetilde{\alpha}_{n;\mathbf{j}}| 
&\leq \frac{1}{1 + Z_{\mathbf{j}}} |\widetilde{\alpha}_{n-1;\mathbf{j}}| + \frac{M_0 \delta t/\epsilon}{1 + Z_{\mathbf{j}}} \nonumber \\
&\leq \left(\frac{1}{1 + Z_{\mathbf{j}}}\right)^n |\widetilde{\alpha}_{0;\mathbf{j}}| + \frac{M_0 \delta t}{\epsilon(1 + Z_{\mathbf{j}})} \sum_{k=0}^{n-1} \left(\frac{1}{1 + Z_{\mathbf{j}}}\right)^k \nonumber \\
&\leq |\widetilde{\alpha}_{0;\mathbf{j}}| + \frac{M_0}{\|\omega_{\mathbf{j}}\|^2+\lambda^2}.
\end{align}

We next bound $\|\widetilde{c}_n\|_{L^\infty}$ by approximating the partial Fourier sum. Using $\|\omega_{\mathbf{j}}\| = \frac{2\pi}{L}\|\mathbf{j}\|$ and treating the sum over $\mathcal{H} \setminus \{\mathbf{0}\}$ as a Riemann sum, which is approximated by an integral in spherical coordinates:
\begin{align}
\|\widetilde{c}_n\|_{L^\infty}
&\leq \frac{1}{L^3}
\sum_{\mathbf{j}\in\mathcal{H}}
|\widetilde{\alpha}_{n;\mathbf{j}}|
\nonumber\\
&\leq \frac{1}{L^3}\left(
\sum_{\mathbf{j}\in\mathcal{H}}
|\widetilde{\alpha}_{0;\mathbf{j}}|
+\frac{M_0}{\lambda^2}
+\sum_{\mathbf{j}\in\mathcal{H}\setminus\{\mathbf{0}\}}
\frac{M_0L^2}{4\pi^2\|\mathbf{j}\|^2}
\right)
\nonumber\\
&\leq \frac{1}{L^3}\left(
\sum_{\mathbf{j}\in\mathcal{H}}
|\widetilde{\alpha}_{0;\mathbf{j}}|
+\frac{M_0}{\lambda^2}
+\frac{M_0L^2}{4\pi^2}
\int_{1}^{\frac{\sqrt{3}}{2}H}
\frac{1}{r^2}\,4\pi r^2\,dr
\right)
\nonumber\\
&\leq C_{\mathrm{stab}}H+C_0,
\end{align}
where $C_{\mathrm{stab}}=\frac{\sqrt{3}M_0}{2\pi L}$ and $C_0=\frac{1}{L^3}\left(
\sum_{\mathbf{j}\in\mathcal{H}}
|\widetilde{\alpha}_{0;\mathbf{j}}|
+\frac{M_0}{\lambda^2}
\right)$ are constants.

Regarding the gradient $\nabla \widetilde{c}$, the SIPF-$r$ method employs a specific discretization to handle the singularity of the Green's function $\mathcal{K}_{\epsilon, \delta t}$. According to Alg.~\ref{alg:SIPF-$r$ho}, the gradient at a particle position $\widetilde{X}^p_n$ is computed as:
\begin{align}\label{nabla c explicit}
\nabla \widetilde{c}(\widetilde{X}^p_n, t_n) &= -\frac{\epsilon}{\delta t} \frac{L^3}{H^3} \sum_{\mathbf{k} \in \mathcal{H}} \nabla_\mathbf{x} \mathcal{K}_{\epsilon,\delta t}(\widetilde{X}^p_n + \bar{X}^p_n - \mathbf{x}_{\mathbf{k}}) \widetilde{c}_{n-1}(\mathbf{x}_{\mathbf{k}} - \bar{X}^p_n) \nonumber \\
&\quad - \frac{M_0}{R}\sum_{s \in C_p, s \neq p} \nabla_{\mathbf{x}} \mathcal{K}_{\epsilon, \delta t}(\widetilde{X}^p_n-\widetilde{X}^s_n),
\end{align}
where the spatial shift $\bar{X}^p_n=\frac{L}{2H}+\lfloor \frac{\widetilde{X}^p_n}{L/H}\rfloor \frac{L}{H}-\widetilde{X}_n^p$ ensures that the evaluation point is bounded away from the singularity of Green's function at grid points $\mathbf{x}_{\mathbf{k}}$. Specifically, let 
$r := \big\| \widetilde{X}^p_n + \bar{X}^p_n - \mathbf{x}_{\mathbf{k}} \big\| \geq \frac{L}{2H}$ and $\beta = \sqrt{\lambda^2 + \epsilon/\delta t}$. The gradient of the Green's function satisfies 
\begin{equation}
\left\| \frac{\epsilon}{\delta t} \nabla \mathcal{K}_{\epsilon,\delta t}(\widetilde{X}^p_n + \bar{X}^p_n - \mathbf{x}_{\mathbf{k}}) \right\| \leq \beta^2 \frac{e^{-\beta r}}{4\pi} \left( \frac{1}{r^2} + \frac{\beta}{r} \right).
\end{equation}
Letting $y = \beta r$, and using the inequality $\sup_{y \geq 0} y^k e^{-y} = (k/e)^k$, we can bound the term to be independent of $\beta$ (and thus independent of $\delta t$):
\begin{align}
\sup_{\beta > 0} \left( \frac{\beta^2}{r^2} e^{-\beta r} + \frac{\beta^3}{r} e^{-\beta r} \right) 
&= \frac{1}{r^4} \sup_{y} (y^2 e^{-y} + y^3 e^{-y}) \nonumber \\
&\leq \frac{1}{(L/2H)^4} \left( \frac{4}{e^2} + \frac{27}{e^3} \right).
\end{align}
Since $\|\widetilde{c}_{n-1}\|_{L^\infty}$ is bounded (as shown above) and the sum over $\mathbf{k} \in \mathcal{H}$ is finite for a fixed $H$, the first term in Eq.~\eqref{nabla c explicit} is uniformly bounded by a constant depending on $H$ but independent of $\delta t$. 

For the second term in Eq.~\eqref{nabla c explicit}, the RBM excludes self-interaction ($s \neq p$). For any fixed grid resolution $H$, the corresponding effective kernel corresponds to a spectrally truncated approximation that is smooth. Thus, for any finite number of particles, this sum is finite.

Combining these results, $\|\nabla \widetilde{c}_n\|_{L^\infty(\Omega)} \leq C_{\mathrm{reg}}(H)$ is guaranteed by the algorithm's design, where the $L^\infty$-norm for the vector-valued gradient is defined by 
$\|\nabla \widetilde{c}_n\|_{L^\infty(\Omega)} := \operatorname{\,sup}_{\mathbf{x}\in\Omega} \|\nabla \widetilde{c}_n(\mathbf{x})\|_{L^2}$, and $C_{\mathrm{reg}}(H)$ is a constant depending on $H$, $L$, $M_0$, and $\lambda$.
\end{proof}

This lemma confirms that the boundedness condition in Assumption \ref{assump:2} is not merely an external hypothesis but a property guaranteed by the SIPF-$r$ method itself. This demonstrates that the SIPF-$r$ method effectively regularizes the singular Keller-Segel kernel. While the exact solution may exhibit finite-time blow-up (where $\|\nabla c\| \to \infty$), the numerical field remains finite for any fixed Fourier mode $H$. This unconditional numerical stability ensures that the algorithm is robust: it enables the simulation of blow-up phenomena by capturing the solution's growth trend as $H$ increases, while avoiding numerical breakdown at fixed resolutions.


\section{Numerical Experiments}\label{section: experiment}
The numerical experiments are organized into three main parts to evaluate the SIPF-$r$ method. In Subsection~\ref{subsection_validation}, we validate the SIPF-$r$ method by quantifying its accuracy against a high-resolution radial finite difference benchmark and verifying its convergence rates with respect to the time step and batch size. Subsection~\ref{subsection:blowup-detection} investigates the method's capability to detect finite-time blow-up phenomena and critical mass thresholds under various conditions. Finally, Subsection~\ref{subsection:assumptions} provides empirical verification of the key theoretical assumptions used in our analysis, specifically the spatial Lipschitz continuity of the concentration gradient.

\subsection{Validation of the SIPF-$r$ Method}\label{subsection_validation}
\subsubsection{Comparison with FDM}\label{subsection: Accuracy of SIPF-$r$ Method} 
We first demonstrate the accuracy of the SIPF-$r$ method. In the radially symmetric case, the fully parabolic KS system \eqref{K-S system} in 3D can be expressed as $\rho(x,y,z,t) = \rho(r,t)$ and $c(x,y,z,t) = c(r,t)$, where $r = \sqrt{x^2 + y^2 + z^2}$. The system is then rewritten as follows in 1D:
\begin{equation}\label{K-S-system-radial}
\left\{
\begin{aligned}
\rho_t &= \mu\, \left (\frac{\partial^2 \rho}{\partial r^2}+\frac{2}{r}\frac{\partial \rho}{\partial r}\right ) - \chi\, \left (\frac{\partial \rho}{\partial r}\frac{\partial c}{\partial r} + \rho\cdot(\frac{\partial^2 c}{\partial r^2}+\frac{2}{r}\frac{\partial c}{\partial r})\right ),  \\
\epsilon c_t &= \, \left (\frac{\partial^2 c}{\partial r^2}+\frac{2}{r}\frac{\partial c}{\partial r}\right ) - \lambda^2 c + \rho.
\end{aligned}
\right.
\end{equation}

The radial representation in 1D allows us to compute a reference solution using a very fine mesh by the finite difference method (FDM). It will serve as a benchmark to quantify the accuracy of the SIPF-$r$ method in 3D. We define the relative error between the cumulative distribution functions (CDFs) obtained from the FDM and the SIPF-$r$ method as
\begin{equation}\label{Relative Error}
\text{Relative Error} = \frac{1}{N} \sum_{i=1}^N 
\begin{cases} 
0, & \text{if } F_{\text{FDM}}(s_i) = 0, \\
\frac{|F_{\text{SIPF-$r$}}(s_i) - F_{\text{FDM}}(s_i)|}{F_{\text{FDM}}(s_i)}, & \text{otherwise,}
\end{cases}
\end{equation}
where $F_{\text{SIPF-$r$}}(s_i)$ and $F_{\text{FDM}}(s_i)$ represent the CDFs of $\rho$ computed via the SIPF-$r$ and FDM methods, respectively, and $s_i$ denotes the $i$-th radial mesh point in the FDM, which are the discrete points along the radial direction starting from the origin. To ensure the relative error is well-defined, we set it to zero wherever $F_{\text{FDM}}(s_i) = 0$.

Here, the initial distribution $\rho_0$ is assumed to be a uniform distribution over a ball centered at $(0,0,0)^T$ with radius 1. The model parameters are chosen as follows:
\begin{align}\label{parameter}
\mu = \chi = 1, \quad \epsilon = 10^{-4}, \quad \lambda = 10^{-1}.
\end{align}
For the numerical computation, we use $H = 24$ Fourier basis functions in each spatial dimension to discretize the chemical concentration $c$ and use $P = 10000$ particles to represent the approximated distribution $\rho$, where the batch size in Algorithm~\ref{alg:SIPF-$r$ho} is $R=\lfloor \sqrt{P} \rfloor=100$. The computational domain is $\Omega = [-L/2, L/2]^3$, where $L = 8$, and the total mass is chosen to be $M_0 = 20$. The evolution of $c$ and $\rho$ is computed using Algorithm~\ref{alg:SIPF} with a time step size $\delta t = 10^{-4}$, up to the final simulation time $T = 0.1$. 

\begin{figure}[htbp]
    \centering
    \begin{subfigure}{0.23\linewidth}
        \centering
\includegraphics[width=\linewidth]{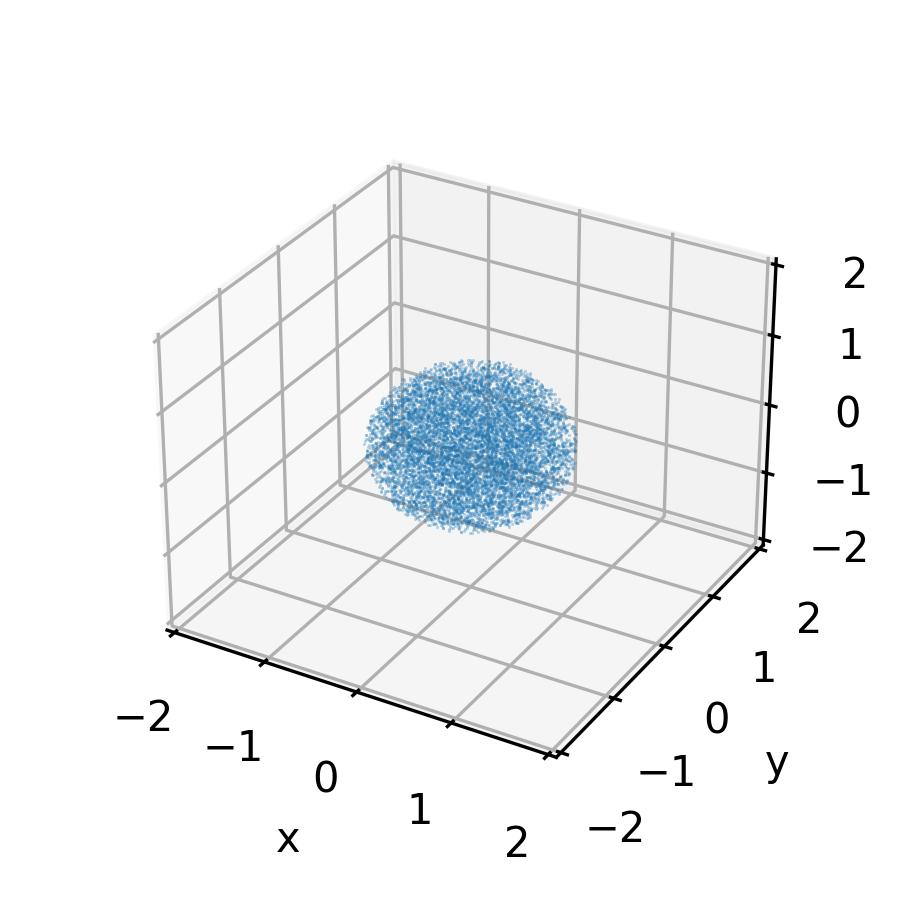}
        \caption{t=0}
    \end{subfigure}%
    ~
    \begin{subfigure}{0.23\linewidth}
        \centering
\includegraphics[width=\linewidth]{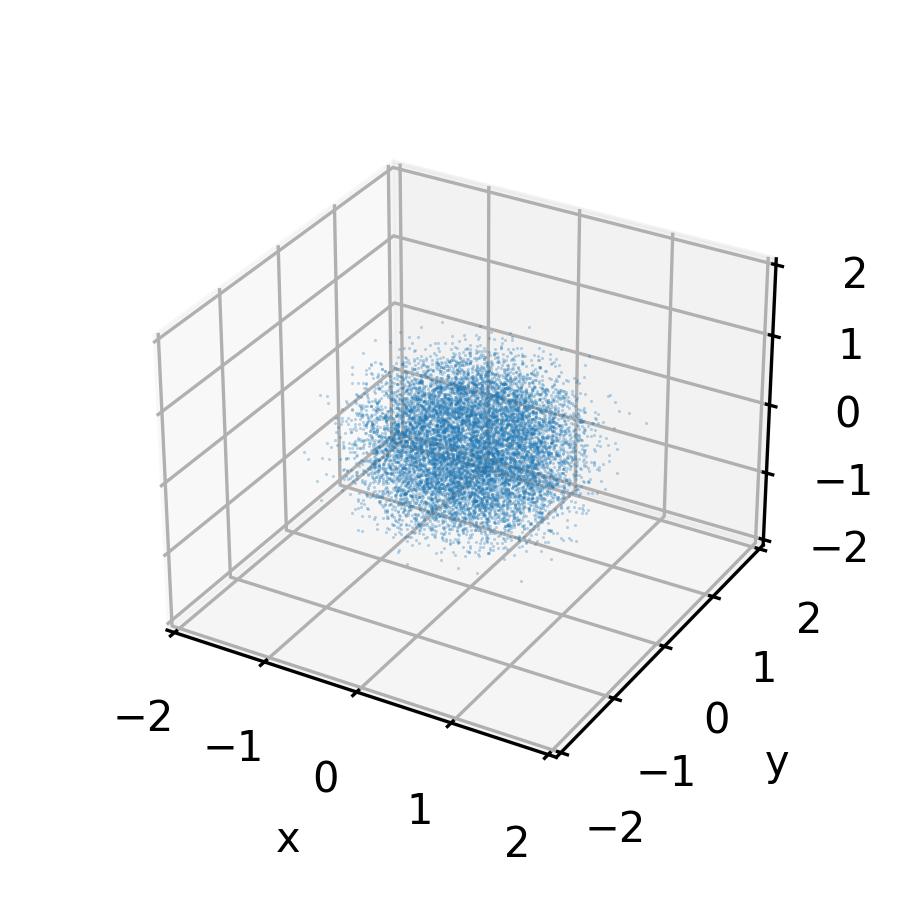}
        \caption{t=0.025}
    \end{subfigure}
    ~
    \begin{subfigure}{0.23\linewidth}
        \centering
\includegraphics[width=\linewidth]{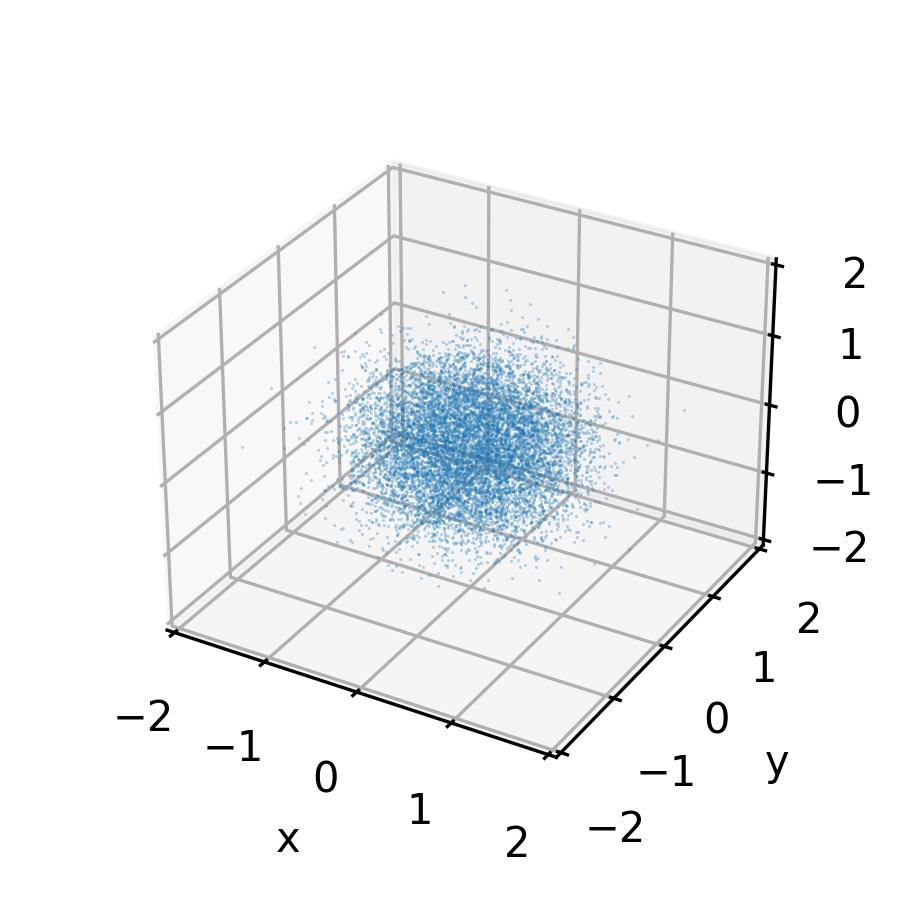}
        \caption{t=0.05}
    \end{subfigure}%
    ~
    \begin{subfigure}{0.23\linewidth}
        \centering
\includegraphics[width=\linewidth]{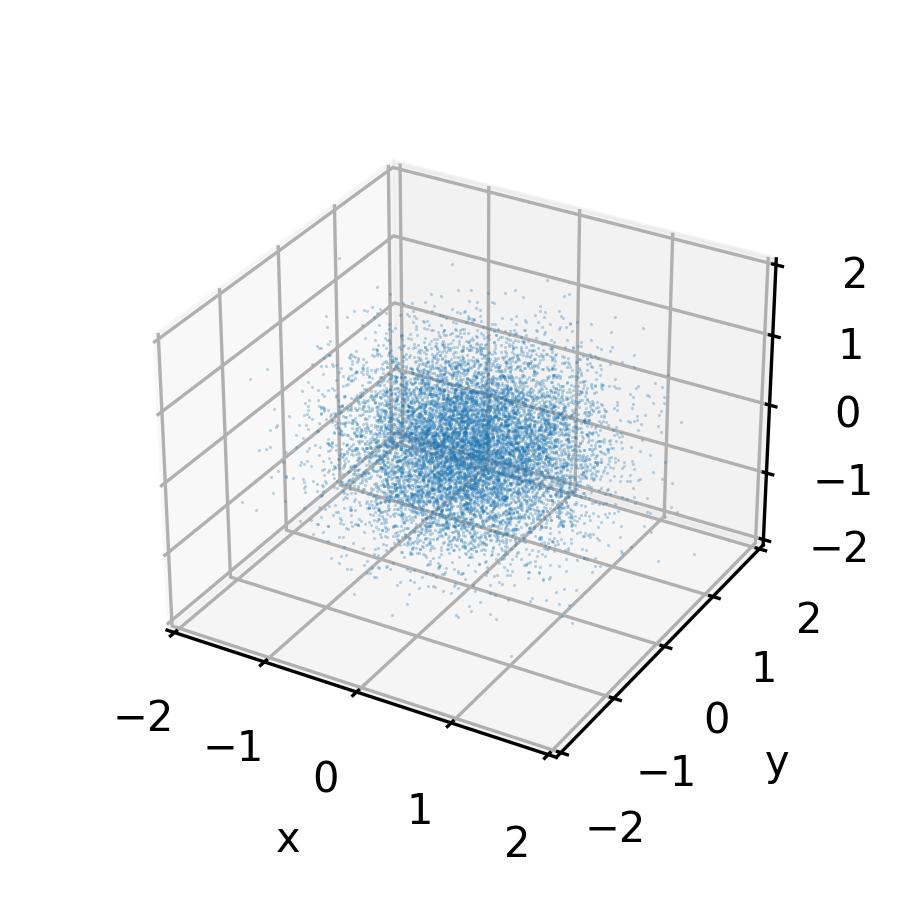}
        \caption{t=0.1}
    \end{subfigure}
    
    \caption{Scatter plot of particles with $M_0$ = 20.}
    \label{evolve}
\end{figure}

\begin{figure}[htbp]  
    \centering  
\includegraphics[width=0.45\textwidth]{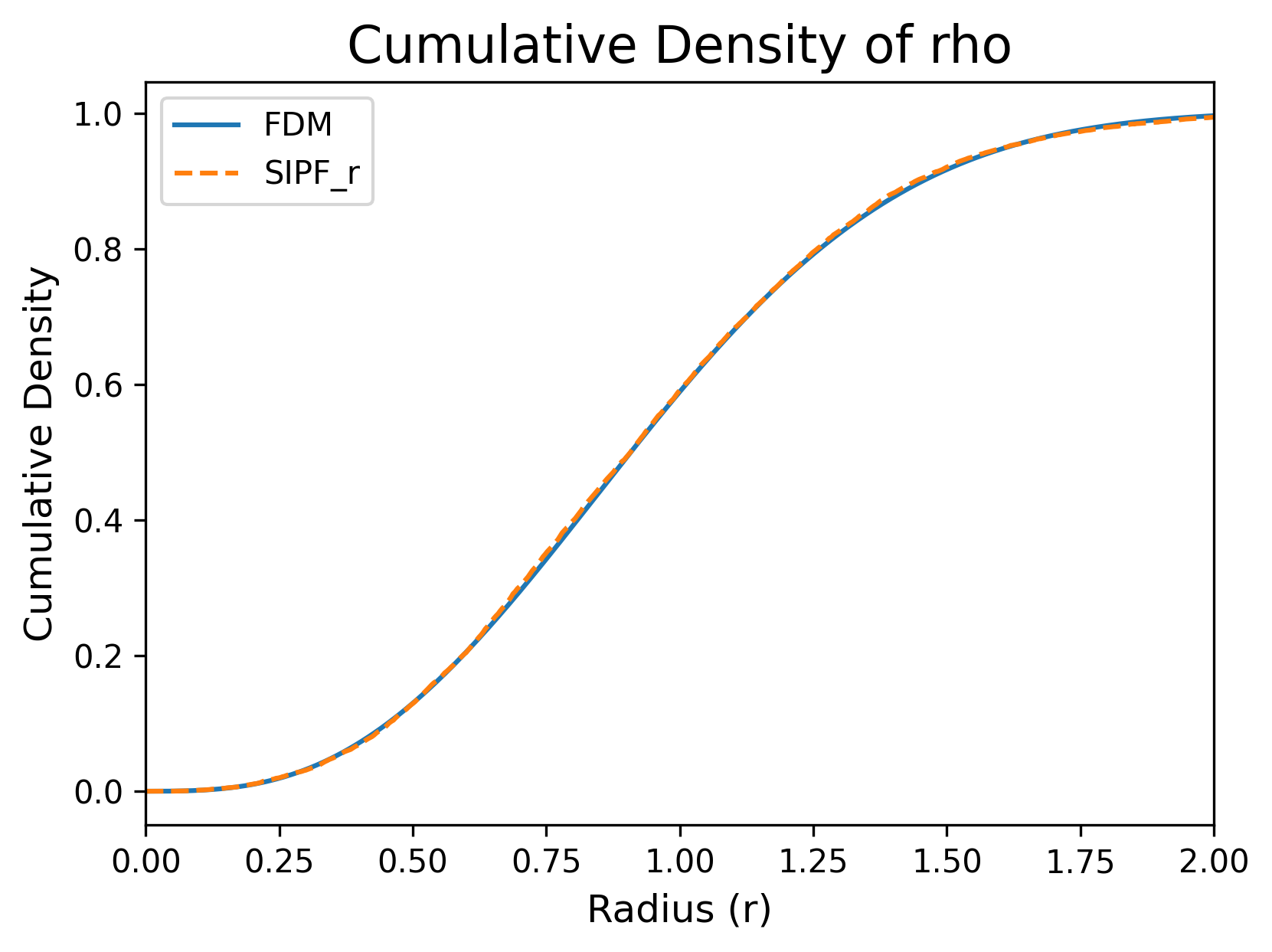}  
    \caption{Cumulative distribution of $\rho$ computed by the SIPF-$r$ method and radial FDM.}  
    \label{comparison}
\end{figure}

In Fig.~\ref{evolve}, we present the evolution of particles over time, showing the dynamic behavior of $\rho$. Additionally, in Fig.~\ref{comparison}, we compare the cumulative probability curves of $\rho$ obtained from the radial FDM and the SIPF-$r$ method at $T = 0.1$, with a mean relative error of 0.05512 as defined in Eq.~\eqref{Relative Error}. This comparison demonstrates that the SIPF-$r$ algorithm achieves high accuracy in approximating the true solution. These results validate the effectiveness of the SIPF-$r$ algorithm in capturing the behavior of the particle distribution.



\subsubsection{Convergence of the SIPF-$r$ Method}\label{convergence_rate}
In this subsection, we validate the convergence of the SIPF-$r$ method numerically. Based on Eq.~\eqref{general alpha}, the error between $\widetilde{c}$ and $c$ can be quantified by the $L^2$ error between their Fourier coefficients $\widetilde{\alpha}$ and $\alpha$. We adopt the same initial conditions as in Subsection \ref{subsection: Accuracy of SIPF-$r$ Method}. To eliminate the uncertainty introduced by the RBM, the reference solution is computed using the original SIPF method \cite{wang2025novel} with parameters $\delta t = 10^{-6}$, $H = 24$, and $P = 10000$. Additionally, we set $M_0=20, T=0.01$ to ensure that the system remains free of singularities, as verified in Fig.~3 of \cite{wang2025novel}.

\begin{figure}[htbp]
    \centering
    \vspace{-0.2cm} 
    
    \begin{subfigure}[b]{0.48\textwidth}
        \centering
        \includegraphics[width=\linewidth]{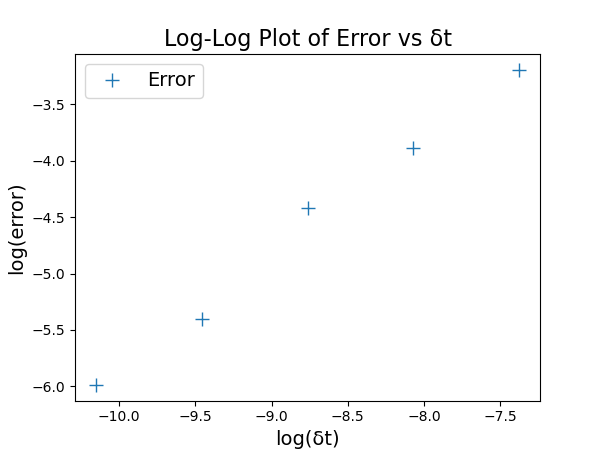}
        \caption{vs. time step $\delta t$ (log-scale)}
        \label{c_error_vs_dt}
    \end{subfigure}
    \hfill
    \begin{subfigure}[b]{0.48\textwidth}
        \centering
        \includegraphics[width=\linewidth]{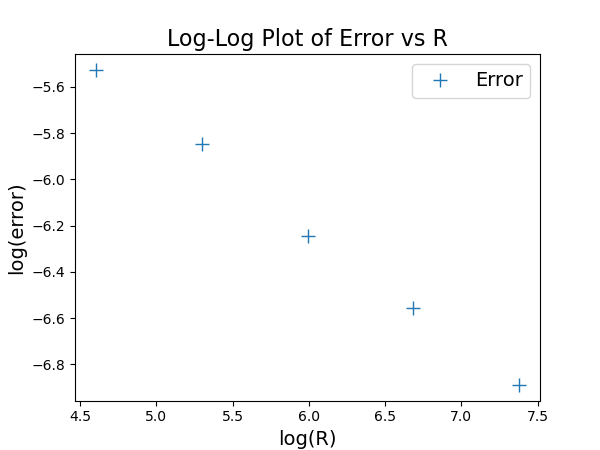}
        \caption{vs. batch size $R$ (log-scale)}
        \label{c_error_vs_R}
    \end{subfigure}
    \caption{$L^2$ error of $\widetilde{c}$ in the SIPF-$r$ method.}
    \label{fig:convergence_rate}
\end{figure}

To investigate the convergence with respect to the time step $\delta t$, we vary $\delta t$ from $2^{-8}T$ to $2^{-4}T$. Since Theorem \ref{theorem} holds with high probability, we perform 100 independent experiments for each $\delta t$ to empirically validate the algorithm's accuracy. The mean $L^2$ error of the Fourier coefficients is computed over these 100 trials. As shown in Fig.~\ref{c_error_vs_dt}, the slope of the mean $L^2$ error versus $\delta t$ on a logarithmic scale indicates an approximate first-order convergence rate, with $e(\delta t) = \mathcal{O}(\delta t^{1.023})$. This result aligns with the theoretical bound given in Eq.~\eqref{theorem inequality} of Theorem \ref{theorem}. 

Furthermore, we examine the mean $L^2$ error of $\widetilde{c}(\cdot, T)$ for varying batch sizes $R = 100, 200, 400, 800, 1600$, while keeping $P=10000$. From Eq.~\eqref{theorem inequality}, with other parameters fixed, the theoretical $L^2$ error of $\widetilde{c}$ with respect to the batch size $R$ should scale as $\mathcal{O}(R^{-\frac{1}{2}})$. This is empirically verified in Fig.~\ref{c_error_vs_R}, where the fitted convergence rate is $e(R) = \mathcal{O}(R^{-0.495})$, closely matching the theoretical prediction.

\subsection{Numerical Detection of Finite-time Blow-up}\label{subsection:blowup-detection}
While the error analysis in Section~\ref{section: proof} assumes regular solutions, our algorithm can effectively detect finite-time blow-up phenomena under practical conditions. The parameters are the same as in Eq.~\eqref{parameter}, with the radial system \eqref{K-S-system-radial} serving as a reference benchmark. We consider a more concentrated initial condition than in the preceding subsections: a uniform distribution within a ball of radius $r = 1.5$.

\begin{figure}[htbp]
\centering
\begin{subfigure}{0.32\textwidth}
\includegraphics[width=\textwidth]{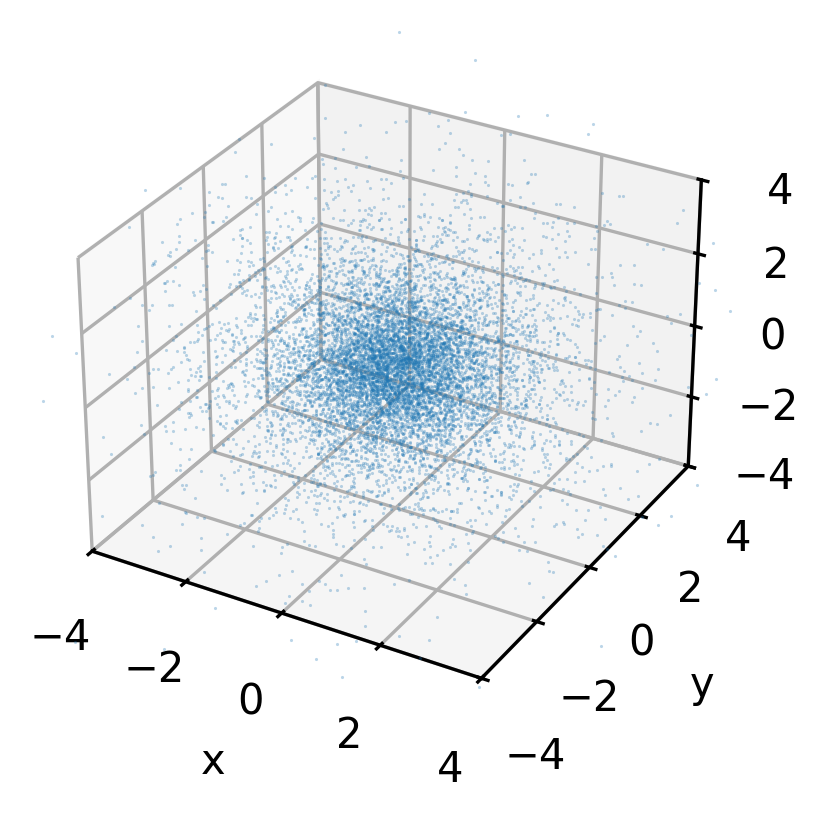}
\caption{Scatter plot of particles at $T=2$ with $M_0=72$.}
\label{fig:blowup-M72}
\end{subfigure}
\hfill
\begin{subfigure}{0.32\textwidth}
\includegraphics[width=\textwidth]{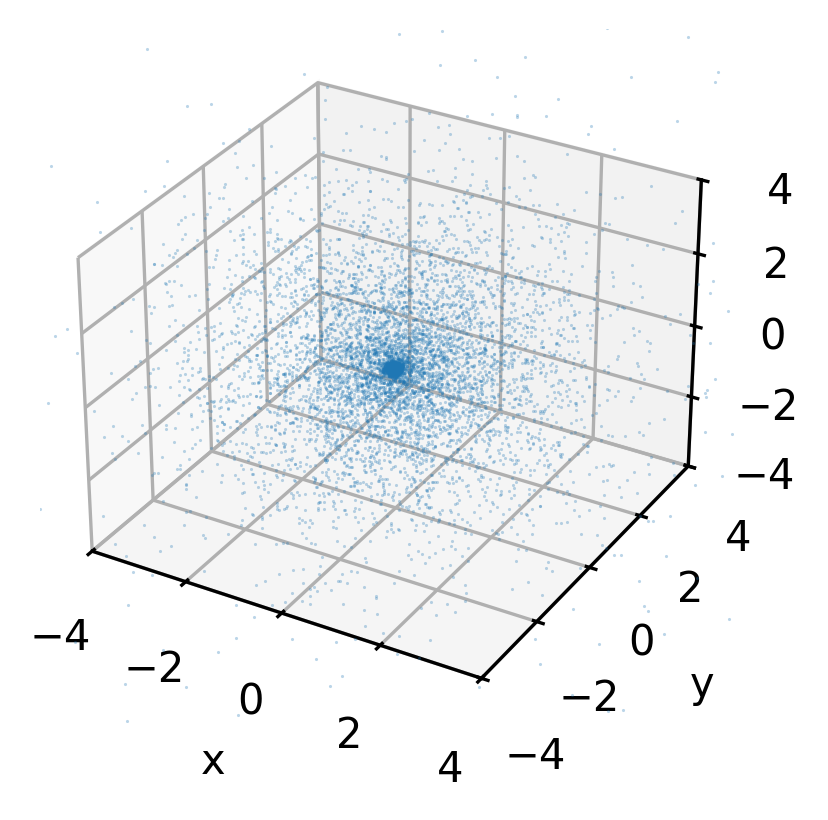}
\caption{Scatter plot of particles at $T=2$ with $M_0=73$.}
\label{fig:blowup-M73}
\end{subfigure}
\hfill
\begin{subfigure}{0.32\textwidth}
\includegraphics[width=\textwidth]{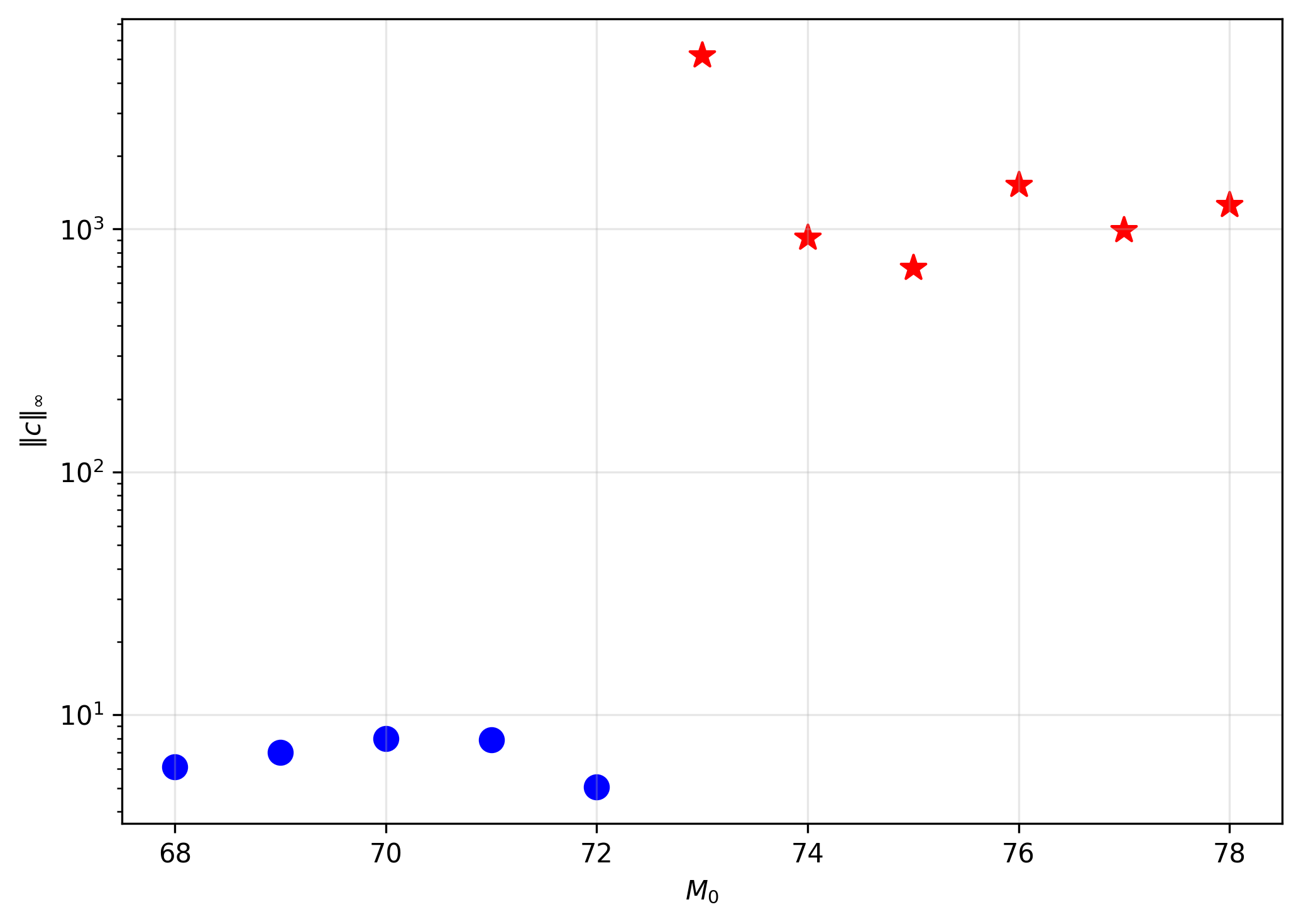}
\caption{$\|c\|_{\infty,\text{FDM}}$ at $T=2$ vs $M_0$ (the red star shape marker denotes the numerical instabilities).}
\label{fig:blowup-FDM}
\end{subfigure}
\caption{Effects of initial mass $M_0$ on concentration behavior (finite time blow-up).}
\label{fig:blowup-comparison}
\end{figure}

Fig.~\ref{fig:blowup-comparison} illustrates the effect of varying the initial mass $M_0$ on solution behavior at $T=2$. \by{We note that unlike the 2D case, the 3D KS system has no universal critical mass, and the blow-up formation depends on both the total mass and the initial spatial distribution. For our prescribed initial configuration, defined as a uniform density supported over a spherical ball of fixed radius, we numerically observed a sharp blow-up transition.} Figs.~\ref{fig:blowup-M72} and \ref{fig:blowup-M73} show particle distributions for $M_0=72$ and $M_0=73$, respectively, with the latter exhibiting clear concentration and potential blow-up. In Fig.~\ref{fig:blowup-FDM}, we present the maximum value of $c$ over space at $T=2$ for different initial masses, computed via the radial FDM as a benchmark. We denote
\begin{equation}\label{eq:c_inf_def}
\|c\|_{\infty,\mathrm{FDM}} = \sup_r |c(r,T)|.
\end{equation}
The FDM exhibits numerical instabilities (marked in red) for initial masses between $72$ and $73$, \by{which identifies the threshold for blow-up specific to this concentrated initial configuration. This agrees with the blow-up transition observed in our proposed method in Fig.~\ref{fig:blowup-M73}.} Remarkably, the concentration and potential blow-up are captured with only $P = 10^4$ particles and a moderate value of $H$, demonstrating the algorithm's ability to detect singularities without the strict constraints imposed in the convergence analysis.

\begin{figure}[htbp]
    \centering
    \begin{subfigure}{0.48\textwidth}
        \centering
        \includegraphics[width=\textwidth]{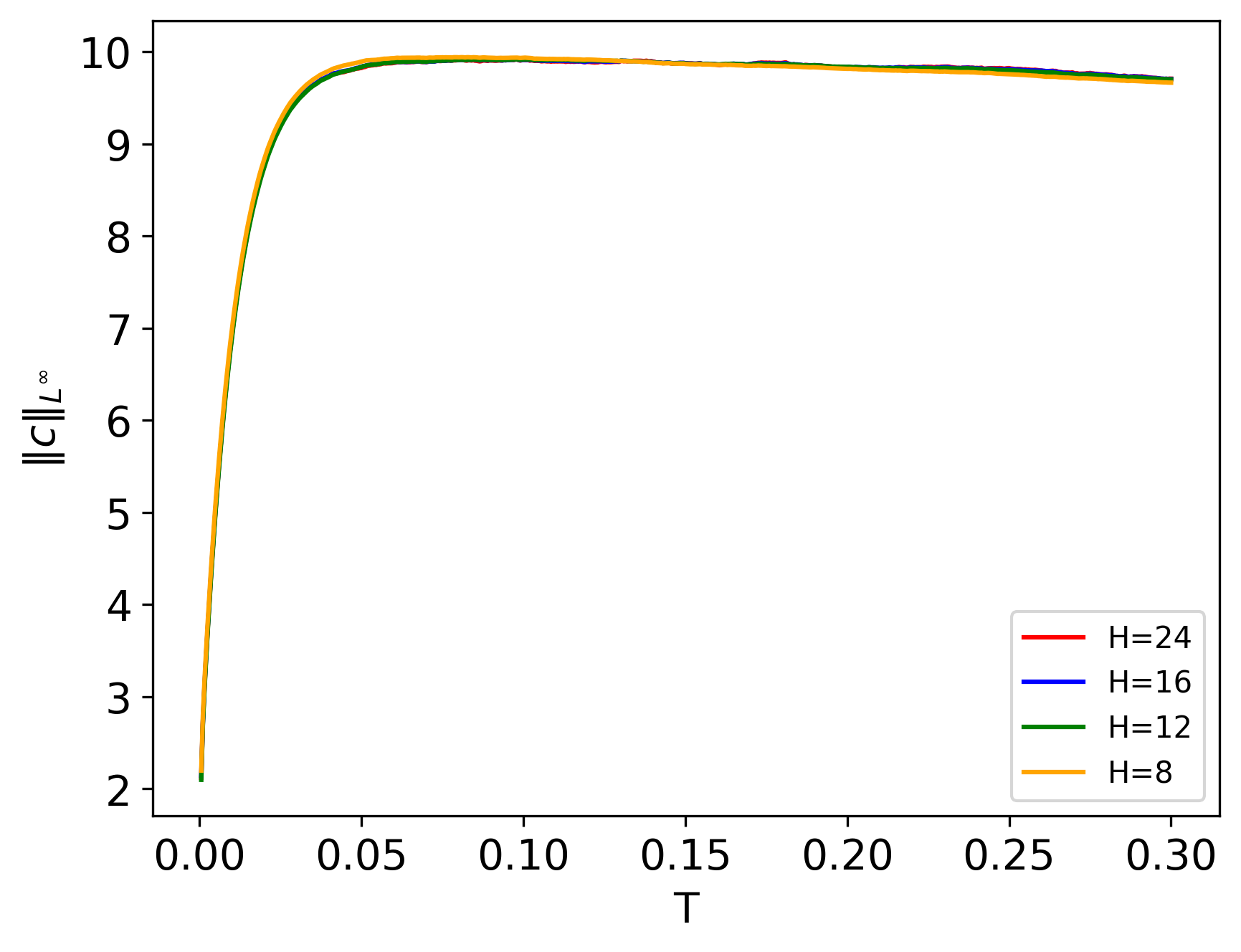}
        \caption{$M_0=40$}
        \label{fig:c_max_M40}
    \end{subfigure}
    \hfill
    \begin{subfigure}{0.48\textwidth}
        \centering
        \includegraphics[width=\textwidth]{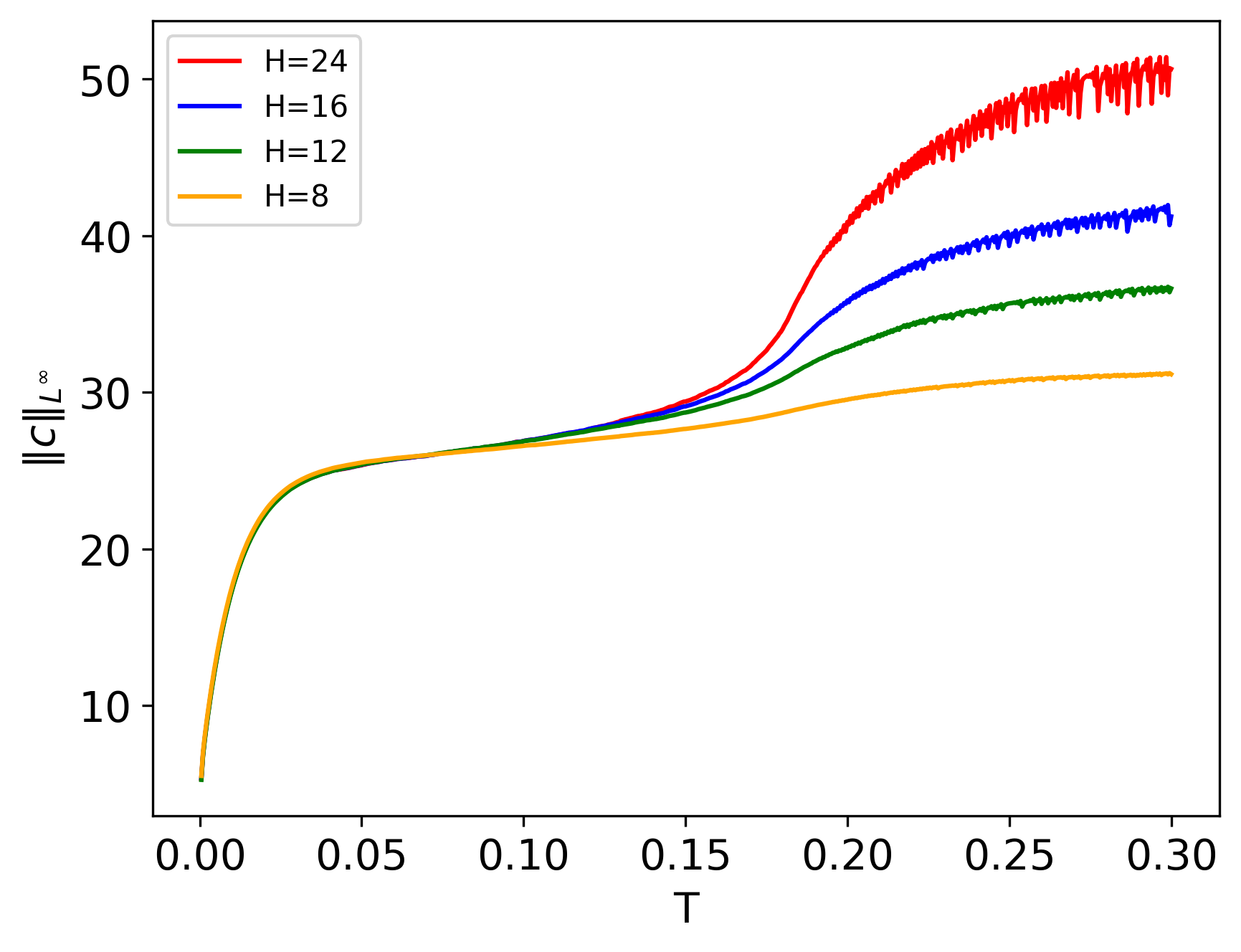}
        \caption{$M_0=100$}
        \label{fig:c_max_M100}
    \end{subfigure}
    \caption{Maximum chemical concentration $c$ vs computation time $T$ for different discretization parameters $H$ and total masses $M_0$}
    \label{fig:c_max_vs_time}
\end{figure}

To further investigate blow-up detection with our method, we examine the maximum chemical concentration $c$ as a function of time $T$ for different discretization levels $H$ and masses $M_0$. As shown in Fig.~\ref{fig:c_max_vs_time}, when $M_0 = 40$ (Fig.~\ref{fig:c_max_M40}), the maximum $c$ remains stable and shows minimal variation across $H$ values, indicating a non-blow-up regime. In contrast, for $M_0 = 100$ (Fig.~\ref{fig:c_max_M100}), the maximum $c$ exhibits strong dependence on $H$, with curves diverging significantly. Notably, this blow-up signature is visible even with coarse discretization ($H=8$ vs $H=12$), demonstrating that singularity detection does not require high-resolution computations. The divergence of solutions at different $H$ levels provides a practical diagnostic for intense concentration and blow-up phenomena without demanding highly refined discretizations.

\begin{figure}[htbp]
\centering
\begin{subfigure}{0.32\textwidth}
\includegraphics[width=\textwidth]{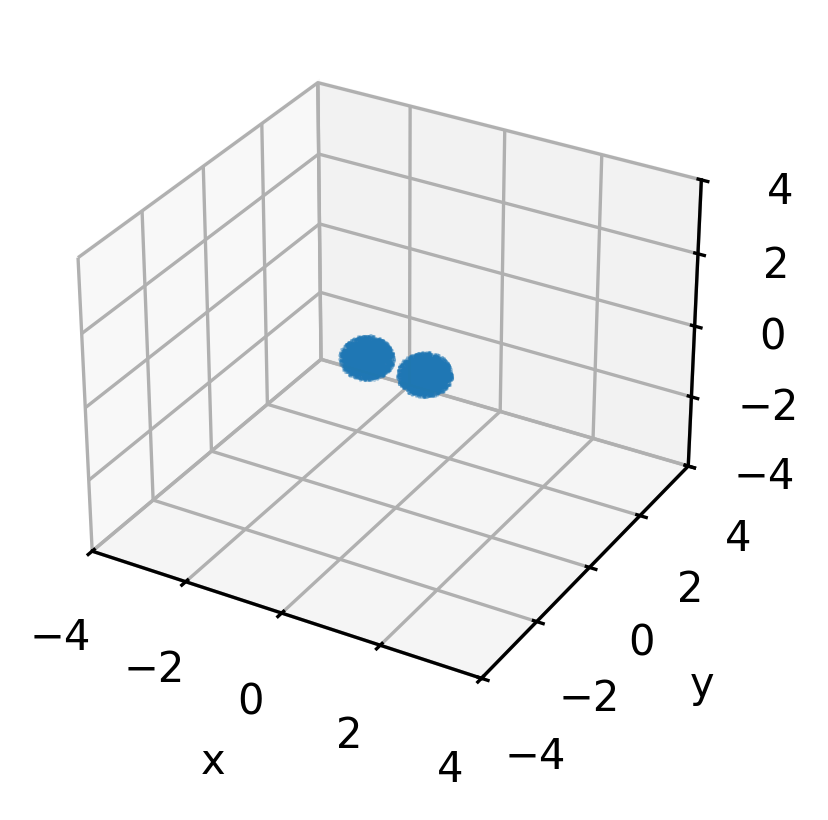}
\caption{$M_0=10$, $T=0$}
\label{fig:two-spheres-M10-T0}
\end{subfigure}
\hfill
\begin{subfigure}{0.32\textwidth}
\includegraphics[width=\textwidth]{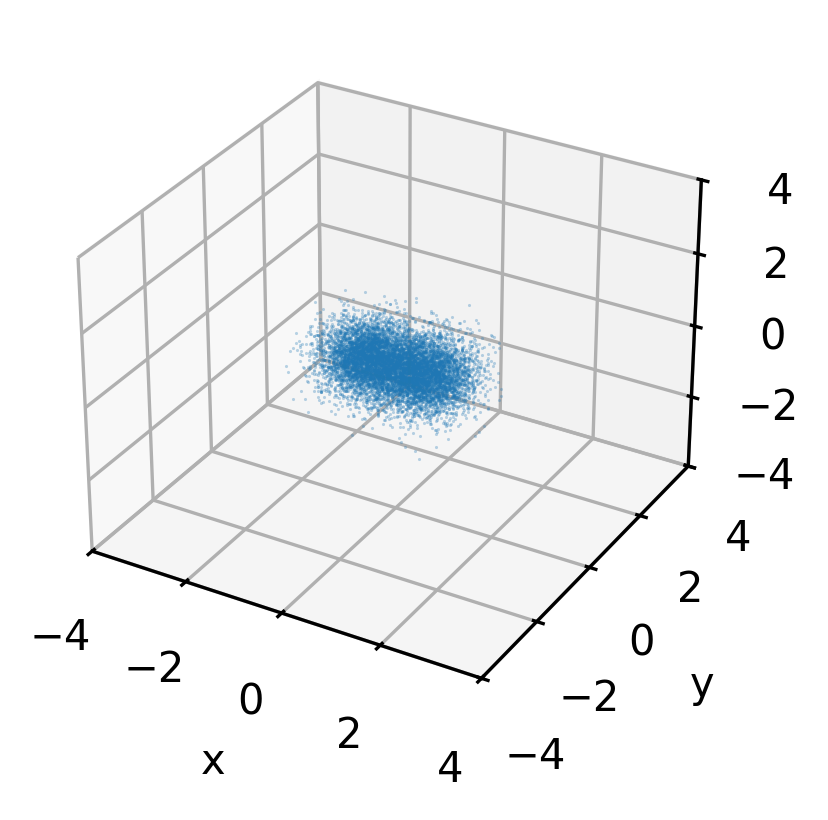}
\caption{$M_0=10$, $T=0.1$}
\label{fig:two-spheres-M10-T1}
\end{subfigure}
\hfill
\begin{subfigure}{0.32\textwidth}
\includegraphics[width=\textwidth]{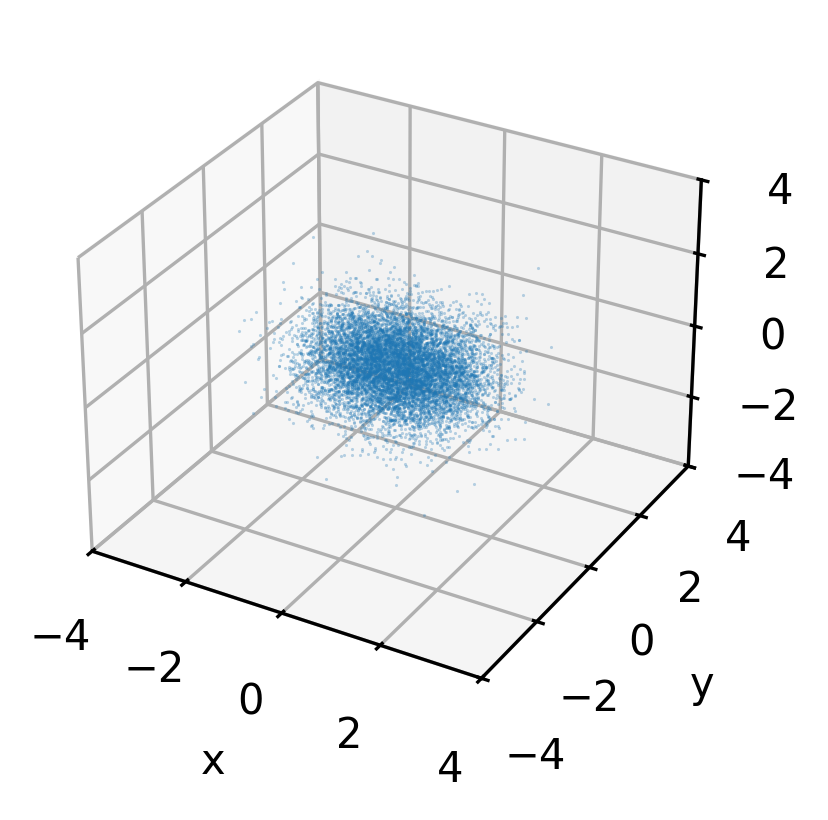}
\caption{$M_0=10$, $T=0.2$}
\label{fig:two-spheres-M10-T2}
\end{subfigure}


\begin{subfigure}{0.32\textwidth}
\includegraphics[width=\textwidth]{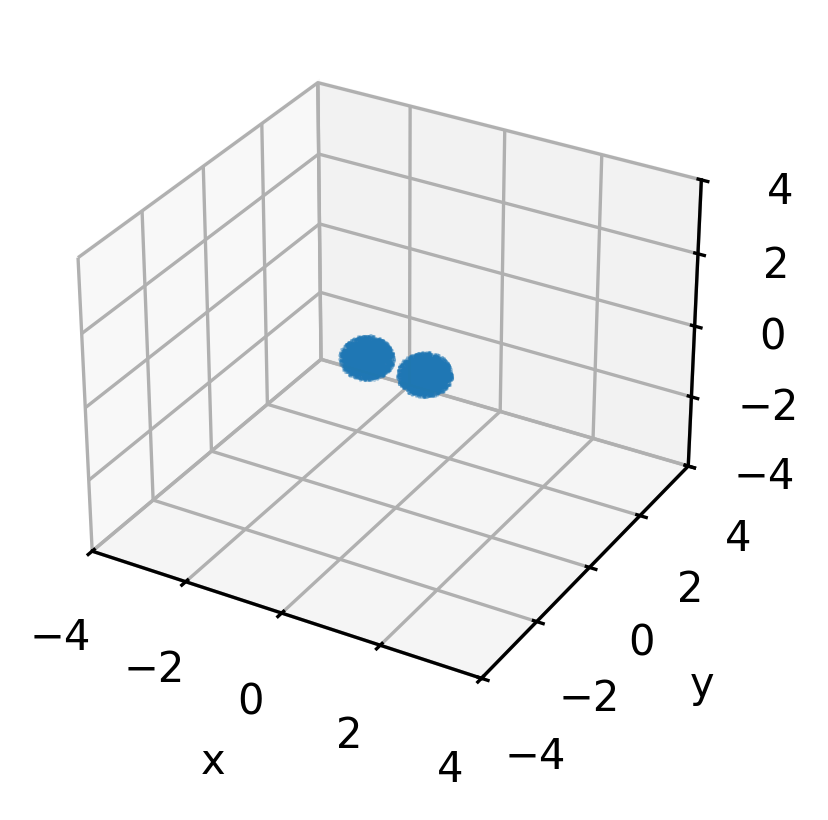}
\caption{$M_0=100$, $T=0$}
\label{fig:two-spheres-M100-T0}
\end{subfigure}
\hfill
\begin{subfigure}{0.32\textwidth}
\includegraphics[width=\textwidth]{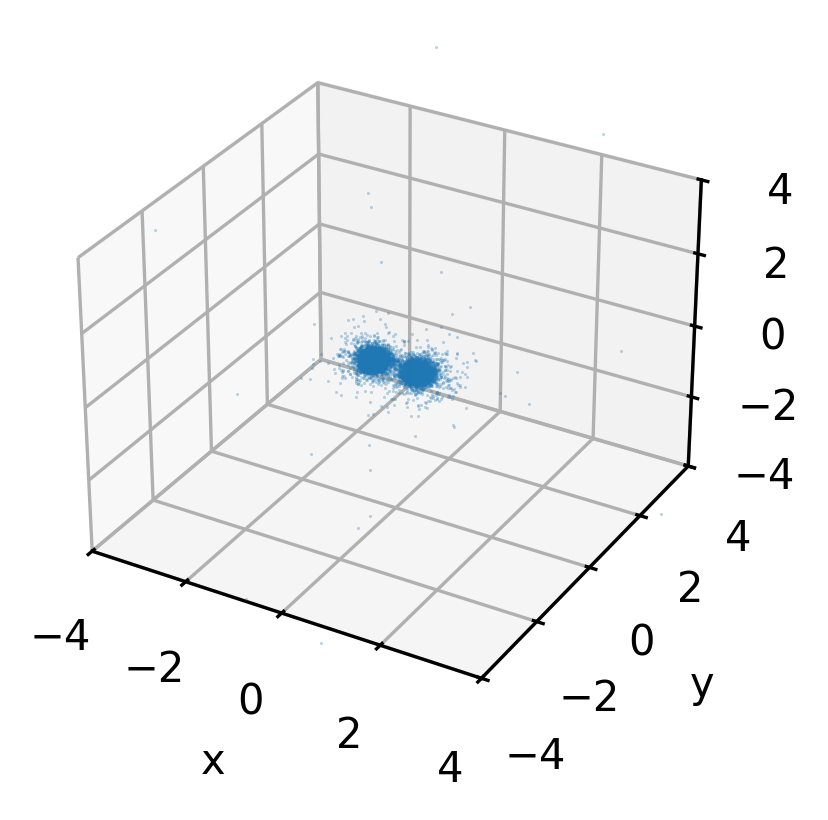}
\caption{$M_0=100$, $T=0.1$}
\label{fig:two-spheres-M100-T1}
\end{subfigure}
\hfill
\begin{subfigure}{0.32\textwidth}
\includegraphics[width=\textwidth]{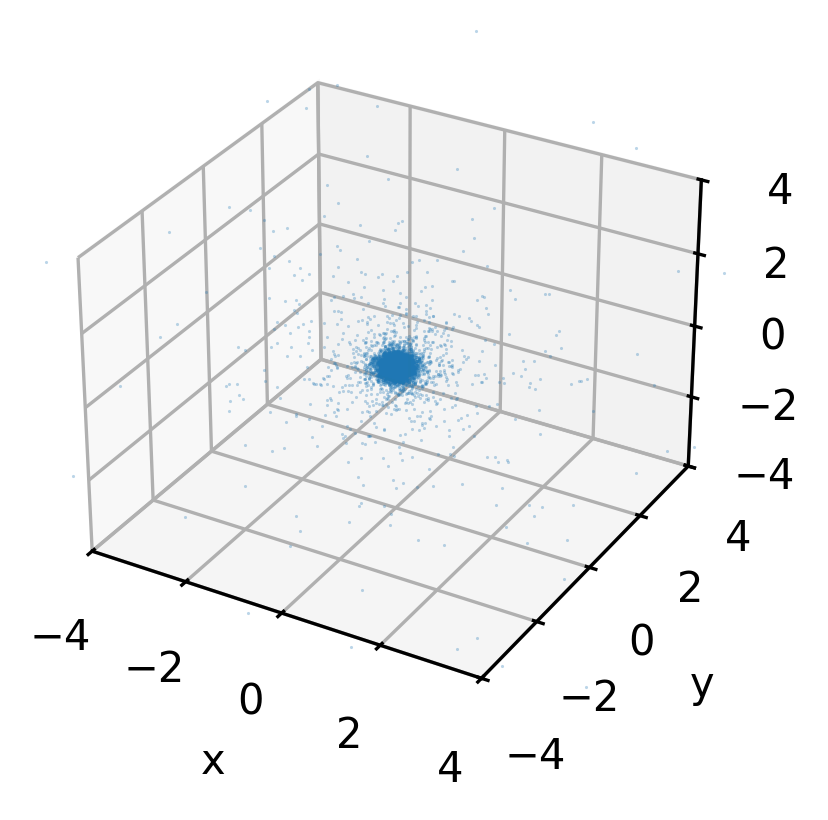}
\caption{$M_0=100$, $T=0.2$}
\label{fig:two-spheres-M100-T2}
\end{subfigure}

\caption{Finite-time blow-up detection for non-radial initial data (two spheres). Top row: $M_0=10$, stable evolution at times $T=0,0.1,0.2$. Bottom row: $M_0=100$, evolution showing blow-up formation at $T=0.2$. Both simulations use $P=10^4$ particles, $H=24$, and initial spheres of radius $r=0.5$ centered at $(\pm0.6,0,0)$.}
\label{fig:two-spheres-blowup}
\end{figure}

We then demonstrate that our method can also capture concentration and potential blow-up for non-radial initial data. We consider two non-overlapping spheres of radius $r=0.5$ with centers at $(0.6, 0, 0)$ and $(-0.6, 0, 0)$, each containing equal mass $M_0/2$ (total mass $M_0=100$). Fig.~\ref{fig:two-spheres-blowup} shows the time evolution for two different mass regimes: the top row ($M_0=10$) remains stable at all times, while the bottom row ($M_0=100$) exhibits clear concentration towards blow-up formation at $T=0.2$. This demonstrates that our algorithm effectively detects blow-up even for non-symmetric initial configurations, highlighting its robustness beyond the radially symmetric case.

\subsection{Validation of Theoretical Assumptions}\label{subsection:assumptions} 
To verify the spatial Lipschitz continuity in Assumption \ref{assump:2}(a), we adjust the spatial discretization by varying $H$ from 6 to 24. At the final time $T=0.1$, we randomly select 1000 pairs of particle points from a total of 10,000 particles in each calculation. The spatial Lipschitz constant $L(H)$ for $\nabla\widetilde{c}$ is defined as the maximum ratio of the gradient difference to the spatial distance over all pairs of particle points $\{\mathbf{x}, \mathbf{y}\}$: 
\begin{equation}\label{L(H)}
L(H) := \max_{\{\mathbf{x}, \mathbf{y}\}} \frac{\|\nabla\widetilde{c}(\mathbf{x}, T) - \nabla\widetilde{c}(\mathbf{y}, T)\|}{\|\mathbf{x} - \mathbf{y}\|}.
\end{equation}
The results, shown in Table \ref{tab:lipschitz_ratios}, list the computed Lipschitz constant $L(H)$ for each value of $H$. The variation in these values is relatively small, confirming that the spatial Lipschitz continuity holds for $\nabla\widetilde{c}$ computed by the SIPF-$r$ method.

\begin{table}[!htbp]
\centering
\begin{tabular}{|c|c|}
\hline
\textbf{Fourier modes ($H$)} & \textbf{Spatial Lipschitz constant ($L(H)$)}\\ \hline
6 & 0.002085 \\ \hline
12 & 0.002106\\  \hline
18 & 0.002036\\ \hline
24 & 0.001957\\ \hline
\end{tabular}
\caption{Spatial Lipschitz constant of $\nabla \widetilde{c}$ vs. $H$.}
\label{tab:lipschitz_ratios}
\end{table}

\byy{
\section{Discussion}\label{section: discussion}

In previous sections, we developed and analyzed the SIPF-$r$ method for the classical 3D fully parabolic KS system \eqref{K-S system}. The underlying SIPF-$r$ framework, however, is flexible and not limited to this classical model. Since the evolution of the chemical field and particle dynamics are updated separately, we can modify the SIPF-$r$ method to incorporate additional biological mechanisms and simulate more complicated mathematical biology models \cite{hu2024stochastic,hu2026novel}. In what follows, we illustrate this flexibility through three representative extensions: multi-species chemotaxis systems with volume-exclusion effects, nonlinear diffusion, and anisotropic cell motility, and we leave their detailed implementation, numerical validation, and convergence analysis for future work.  

\textit{Multi-species chemotaxis system with volume-exclusion effects.}
The SIPF-$r$ framework is readily extendable to multi-species chemotaxis with volume-exclusion (often referred to as volume-filling) effects, which are essential for modeling realistic biological scenarios where cells of different species compete for finite physical space \cite{hillen2001global,painter2002volume,burger2006keller}. A representative two-species chemotaxis model with cross-volume exclusion is 
\begin{equation}\label{eq:discussion_two_species_ve}
\begin{aligned}
\partial_t \rho_k
&=\nabla\cdot\left(D_k(\rho_{\mathrm{tot}})\nabla\rho_k
-\chi_k\rho_k\Phi(\rho_{\mathrm{tot}})\nabla c_k\right),\\
\epsilon\partial_t c_k
&=\Delta c_k-\lambda_k^2c_k+\sum_{l=1}^2\beta_{kl}\rho_l,
\qquad k=1,2,
\end{aligned}
\end{equation}
where $\rho_{\mathrm{tot}} = \rho_1 + \rho_2$ denotes the total cell density, $c_k$ is the chemical concentration of species $k$, $\epsilon$, $\chi_k$, and $\lambda_k$ are positive constants, $\beta_{kl}$ are the production rates of chemical $k$ by species $l$, $D_k(\rho_{\mathrm{tot}})$ is the density-dependent diffusivity, and $\Phi(\rho_{\mathrm{tot}}) = (1 - \rho_{\mathrm{tot}}/\rho_{\max})_+$ is the crowding factor (derived from the volume-filling probability in \cite{painter2002volume}) that modulates chemotactic mobility. This formulation naturally extends to $m$ species by redefining $\rho_{\mathrm{tot}} = \sum_{l=1}^m \rho_l$. Crucially, the system ~\eqref{eq:discussion_two_species_ve} reduces to the classical multi-species KS model when $D_k \equiv \mu_k$ and $\Phi \equiv 1$.

Since the equations governing the chemical field $c_k$ remain linear, we apply Algorithm~\ref{alg:SIPF-C} to update the Fourier coefficients of the chemical field $c_k$ as follows:
\begin{equation}\label{eq:discussion_ve_fourier_update}
\widetilde{\alpha}_{n;\mathbf{j}}^k
=\frac{\frac{\epsilon}{\delta t}\widetilde{\alpha}_{n-1;\mathbf{j}}^k
+\sum_{l=1}^m\beta_{kl}\widehat{\widetilde{\rho}}_{l,n;\mathbf{j}}}
{\epsilon/\delta t+\|\omega_{\mathbf{j}}\|^2+\lambda_k^2},
\qquad k=1,\ldots,m,
\end{equation}
where $\widetilde{\alpha}_{n;\mathbf{j}}^k$ denotes the Fourier coefficient of the chemical field $c_k$, $\delta t$ is the time step, $\omega_{\mathbf{j}} = \frac{2\pi}{L}\mathbf{j}$ is the frequency vector associated with the multi-index $\mathbf{j}=(j_1,j_2,j_3)$ (with $L$ denoting the domain size), $\mathbf{j}\in \mathcal{H}$ (the index set defined in \eqref{index set H}), and $\widehat{\widetilde{\rho}}_{l,n;\mathbf{j}}$ is the Fourier coefficient of the empirical density of species $l$ with $1\leq l\leq m$.

For the density update, species $l$ with total mass $M_l$ is represented at time $t_n$ using $P_l$ particles $\{\widetilde{X}_{l,n}^p\}_{p=1}^{P_l}$, with the corresponding empirical measure
\begin{equation}\label{eq:discussion_ve_empirical_density}
\widetilde{\rho}_{l,n}(\mathbf{x})
=\frac{M_l}{P_l}\sum_{p=1}^{P_l}
\delta(\mathbf{x}-\widetilde{X}_{l,n}^p).
\end{equation}

The Fourier coefficients of the density source in Eq.~\eqref{eq:discussion_ve_fourier_update} can be computed directly from this empirical measure as
\begin{equation}
\widehat{\widetilde{\rho}}_{l,n;\mathbf{j}}
=\mathcal F_{\mathbf{j}}[\widetilde{\rho}_{l,n}]
=\frac{M_l}{P_l}\sum_{p=1}^{P_l}
e^{-i\omega_{\mathbf{j}}\cdot\widetilde{X}_{l,n}^p}.
\end{equation}

The density of species $l$ is then reconstructed at the Fourier grid points via the inverse Fourier transform, and the total density is obtained by summing these reconstructed densities over all species:
\begin{equation}\label{eq:discussion_ve_density_reconstruction}
\widetilde{\rho}_{l,n}^{H}(\mathbf{x})
=\frac{1}{L^3}\sum_{\mathbf{j}\in\mathcal H}
\widehat{\widetilde{\rho}}_{l,n;\mathbf{j}}
e^{i\omega_{\mathbf{j}}\cdot\mathbf{x}},
\qquad
\widetilde{\rho}_{\mathrm{tot},n}^{H}(\mathbf{x})
=\sum_{l=1}^m\widetilde{\rho}_{l,n}^{H}(\mathbf{x}).
\end{equation}

Using the reconstructed total density, the functions
$D_{k,n}(\mathbf{x})=D_k(\widetilde{\rho}_{\mathrm{tot},n}^{H}(\mathbf{x}))$ and
$\Phi_n(\mathbf{x})=\Phi(\widetilde{\rho}_{\mathrm{tot},n}^{H}(\mathbf{x}))$
are first evaluated at the Fourier grid points $x_{\mathbf q}=L\mathbf q/H$, $\mathbf q\in\mathcal H$. The unnormalized Fourier coefficients of these functions are approximated as
\begin{equation}\label{eq:discussion_ve_field_coefficients}
\begin{aligned}
\widehat{D}_{k,n;\mathbf j}
&=\frac{L^3}{H^3}\sum_{\mathbf q\in\mathcal H}
D_{k,n}(x_{\mathbf q})e^{-i\omega_{\mathbf j}\cdot x_{\mathbf q}},\\
\widehat{\Phi}_{n;\mathbf j}
&=\frac{L^3}{H^3}\sum_{\mathbf q\in\mathcal H}
\Phi_n(x_{\mathbf q})e^{-i\omega_{\mathbf j}\cdot x_{\mathbf q}},
\qquad \mathbf j\in\mathcal H.
\end{aligned}
\end{equation}
This coefficient computation extends Step~5 of Algorithm~\ref{alg:SIPF-C} to handle the density-dependent fields $D_{k,n}$ and $\Phi_n$. The values required for the particle update are obtained via Fourier interpolation:
\begin{equation}\label{eq:discussion_ve_coeff_interp}
\begin{aligned}
D_{k,n}(\widetilde{X}_{k,n}^p)
&=\frac{1}{L^3}\sum_{\mathbf{j}\in\mathcal H}\widehat{D}_{k,n;\mathbf{j}}
e^{i\omega_{\mathbf{j}}\cdot\widetilde{X}_{k,n}^p},\\
\nabla D_{k,n}(\widetilde{X}_{k,n}^p)
&=\frac{1}{L^3}\sum_{\mathbf{j}\in\mathcal H}(i\omega_{\mathbf{j}})
\widehat{D}_{k,n;\mathbf{j}}e^{i\omega_{\mathbf{j}}\cdot\widetilde{X}_{k,n}^p},\\
\Phi_n(\widetilde{X}_{k,n}^p)
&=\frac{1}{L^3}\sum_{\mathbf{j}\in\mathcal H}\widehat{\Phi}_{n;\mathbf{j}}
e^{i\omega_{\mathbf{j}}\cdot\widetilde{X}_{k,n}^p}.
\end{aligned}
\end{equation}

These quantities, determined by the density at time $t_n$, are kept fixed during the update from $t_n$ to $t_{n+1}$. Assuming that the frozen diffusivity $D_{k,n}$ is sufficiently smooth and nonnegative (after regularization if needed), the Fokker--Planck correspondence leads to the following update for the $p$-th particle of species $k$:
\begin{equation}\label{eq:discussion_ve_particle_update}
\begin{aligned}
\widetilde{X}_{k,n+1}^p={}&\widetilde{X}_{k,n}^p
+\nabla D_{k,n}(\widetilde{X}_{k,n}^p)\delta t
+\sqrt{2D_{k,n}(\widetilde{X}_{k,n}^p)\delta t}\,N_{k,n}^p
-\chi_k\Phi_n(\widetilde{X}_{k,n}^p)\Biggl[\\
&\quad\nabla\mathcal K_{\epsilon,\delta t}^k
\ast\left(\frac{\epsilon}{\delta t}\widetilde{c}_{k,n-1}\right)(\widetilde{X}_{k,n}^p)
+\sum_{l=1}^m\frac{\beta_{kl}M_l}{R_l}\sum_{s\in C_{p,l}}
\nabla\mathcal K_{\epsilon,\delta t}^k
(\widetilde{X}_{k,n}^p-\widetilde{X}_{l,n}^s)
\Biggr]\delta t,
\end{aligned}
\end{equation}
where $N_{k,n}^p\sim\mathcal N(\mathbf{0},I_3)$ denote independent standard Gaussian vectors, and $\mathcal K_{\epsilon,\delta t}^k$ is the Green's function for the operator $\Delta-\lambda_k^2-\epsilon/\delta t$. For each $l$, the batch $C_{p,l}$ consists of $R_l$ indices sampled with replacement from $\{1,\ldots,P_l\}$; when $l=k$, any occurrence with $s=p$ is omitted from the sum. The bracketed expression approximates $-\nabla\widetilde c_{k,n}$, with its density contribution evaluated by random batching. Using Eq.~\eqref{eq:discussion_ve_empirical_density}, the updated positions yield the empirical densities and hence the source coefficients for the subsequent chemical field update. If $D_k$ is constant and $\Phi\equiv 1$, Eq.~\eqref{eq:discussion_ve_particle_update} reduces to the standard SIPF-$r$ update.

\textit{Chemotaxis systems with nonlinear diffusion.} 
Density-dependent diffusion is widely used in chemotaxis models to describe changes in cellular dispersal induced by the local population density, particularly when crowding leads to nonlinear motility effects \cite{kowalczyk2005preventing,eberl2001new,hillen2009user,kuiper2007global,hofer1995dictyostelium}. For the porous-medium diffusion term $\Delta(\rho^\gamma)=\nabla\cdot(\gamma\rho^{\gamma-1}\nabla\rho)$ with $\gamma>1$, the diffusivity at time $t_n$ can be approximated by $D_n=\gamma\widetilde{\rho}_n^{\gamma-1}$. The values of $D_n$ and $\nabla D_n$ evaluated at particle positions are then computed from the density reconstruction in Eq.~\eqref{eq:discussion_ve_density_reconstruction} and the Fourier interpolation in Eq.~\eqref{eq:discussion_ve_coeff_interp}. After freezing the regularized $D_n$ over one substep, we formally represent the diffusion equation using the associated It\^o SDE, discretize it via the Euler--Maruyama scheme, and then apply the original chemotactic update. 

\textit{Chemotaxis systems with anisotropic cell motility.}
Anisotropic motility arises in structured environments such as extracellular-matrix fibers and neural tracts \cite{othmer2002diffusion,painter2009continuous}. It can be modeled as
\begin{equation}\label{eq:discussion_anisotropic_rho}
\partial_t\rho
=\nabla\cdot\left(\mathbb{D}(\mathbf{x})\nabla\rho
-\rho\chi(\mathbf{x})\nabla c\right),
\end{equation}
where $\mathbb{D}(\mathbf{x})=(D_{ij}(\mathbf{x}))_{ij}$ denotes a sufficiently smooth, symmetric positive semidefinite motility tensor. If chemical diffusion remains isotropic, the spectral update of $\widetilde c$ remains unchanged. The Fokker--Planck correspondence yields the formal particle update
\begin{equation}\label{eq:discussion_anisotropic_particle}
\widetilde{X}_{n+1}^p=\widetilde{X}_n^p
+\left[(\nabla\cdot\mathbb{D})(\widetilde{X}_n^p)
+\chi(\widetilde{X}_n^p)\nabla\widetilde{c}(\widetilde{X}_n^p,t_n)\right]\delta t+\sqrt{2\mathbb{D}(\widetilde{X}_n^p)\delta t}\,N_n^p,
\end{equation}
where $[(\nabla\cdot\mathbb{D})(\mathbf{x})]_i=\sum_j\partial_{x_j}D_{ij}(\mathbf{x})$, $\sqrt{2\mathbb{D}}$ denotes the matrix square root,
and $N_{n}^p\sim\mathcal N(\mathbf{0},I_3)$ denote independent standard Gaussian vectors. Thus, anisotropy modifies the particle solver, whereas the field update and random-batch evaluation of $\nabla\widetilde c$ preserve the structure of Algorithm~\ref{alg:SIPF-$r$ho}. We leave the interpolation of $\mathbb D$ and its divergence, together with the convergence analysis of the resulting scheme, for future work in a separate implementation study.
}

\section{Conclusions}\label{section: conclusion}
In this paper, we introduced a random batch variant \cite{jin2020random} of the original SIPF method \cite{wang2025novel} to simulate the 3D fully parabolic KS system. This modification leverages the randomness in batch sampling to bypass the mean-field limit, which reduces computational complexity without sacrificing accuracy. We established a comprehensive convergence analysis for the fully coupled particle-field system. Specifically, we proved the convergence with high probability for both the density $\widetilde{\rho}(\mathbf{x}, t)$ and the concentration field $\widetilde{c}(\mathbf{x}, t)$ to their respective exact solutions ${\rho}(\mathbf{x}, t)$ and $c(\mathbf{x}, t)$. The error bounds revealed a dependence on $\delta t$, $H$, and $P$, with the density and concentration fields exhibiting distinct but interrelated convergence behaviors. 

Finally, we performed numerical experiments to verify these theoretical rates and confirm the robustness of the SIPF-$r$ method. \by{Notably, the algorithm effectively captures density concentration and potential finite-time blow-up phenomena subject to the critical threshold of initial mass in three space dimensions under specific initial data configurations, even with discretization parameters milder than those required by the worst-case theoretical bounds}. 

Future work will focus on refining the numerical analysis of the classical model and extending the SIPF-$r$ method to broader biological applications. From a theoretical perspective, we will refine the error estimates, particularly to weaken the strong dependence on the Fourier mode $H$, and improve the efficiency of the algorithm in high-dimensional settings. \byy{Beyond refinements to the classical model, we will build on the formulations presented in Section~\ref{section: discussion} to develop SIPF-$r$ schemes for multi-species chemotaxis systems with volume-exclusion effects, nonlinear diffusion, and anisotropic cell motility, and to investigate their numerical performance and convergence properties.}


\section*{Acknowledgements}
\noindent ZZ was partially supported by the National Natural Science Foundation of China (Projects 92470103 and 12171406), the Hong Kong RGC grant (projects 17304324 and 17300325), the Seed Funding Programme for Basic Research (HKU), and the Hong Kong RGC Research Fellow Scheme 2025.  ZW was partially supported by NTU SUG-023162-00001 and MOE AcRF Tier 1 Grant RG17/24. JX was partially supported by NSF grant DMS-2309520, the Swedish Research Council grant no. 2021-06594 at the Institut Mittag-Leffler in Djursholm, Sweden, and the E. Schr\"odinger Institute in Vienna, Austria, during his stay in the Fall of 2025.  The simulations were performed on the research computing facilities of the Information Technology Services at the University of Hong Kong, and the Greenplanet Cluster at UC Irvine.

\bibliographystyle{siamplain}
\bibliography{reference}
\end{document}

%% file: ex_shared.tex
\usepackage{lipsum}
\usepackage{amsfonts}
\usepackage{graphicx}
\usepackage{epstopdf}
\usepackage{algorithmic}
\ifpdf
  \DeclareGraphicsExtensions{.eps,.pdf,.png,.jpg}
\else
  \DeclareGraphicsExtensions{.eps}
\fi

\newsiamremark{remark}{Remark}
\newsiamremark{hypothesis}{Hypothesis}
\crefname{hypothesis}{Hypothesis}{Hypotheses}
\newsiamthm{claim}{Claim}
\newsiamremark{fact}{Fact}
\crefname{fact}{Fact}{Facts}

\headers{Convergence Analysis of SIPF algorithm for Keller-Segel Systems}{Boyi Hu, Zhongjian Wang, Jack Xin, Zhiwen Zhang}

\title{Convergence Analysis of a Stochastic Interacting Particle-Field Algorithm for 3D Parabolic-Parabolic Keller-Segel Systems\thanks{To appear in Mathematics of Computation.
}}

\author{Boyi Hu\thanks{Department of Mathematics, The University of Hong Kong, Pokfulam Road, Hong Kong SAR, P.R.China. huby22@connect.hku.hk.}
\and 
Zhongjian Wang \thanks{Division of Mathematical Sciences, School of Physical and Mathematical Sciences, Nanyang Technological University, 21 Nanyang Link, Singapore 637371. zhongjian.wang@ntu.edu.sg}
\and
Jack Xin \thanks{Corresponding author. Department of Mathematics, University of California at Irvine, Irvine, CA 92697, USA. jack.xin@uci.edu}
\and
Zhiwen Zhang\thanks{Corresponding author. Department of Mathematics, The University of Hong Kong, Pokfulam Road, Hong Kong SAR, P.R.China. Materials Innovation Institute for Life Sciences and Energy (MILES), HKU-SIRI, Shenzhen, 518045, P.R. China. zhangzw@hku.hk.}
}

\usepackage{amsopn}
